\documentclass{article}
\usepackage{amsmath, amssymb}
\usepackage[all]{xy}

\newtheorem{thm}{Theorem}[section]
\newtheorem{cor}[thm]{Corollary}
\newtheorem{prop}[thm]{Proposition}
\newtheorem{lemma}[thm]{Lemma}
\newenvironment{dfn}{\medskip\refstepcounter{thm}
\noindent{\bf Definition \thesection.\arabic{thm}\ }}{\medskip}
\newenvironment{proof}[1][,]{\medskip\ifcat,#1
\noindent{{\it Proof}:\ }\else\noindent{\it Proof of #1.\ }\fi}
{\hfill$\square$\medskip}
\newenvironment{remark}[1][Remark]{\begin{trivlist}
\item[\hskip \labelsep {\bfseries #1}]}{\end{trivlist}}

\newenvironment{note}[1][Note]{\begin{trivlist}
\item[\hskip \labelsep {\bfseries #1}]}{\end{trivlist}}
\newenvironment{notes}[1][Notes]{\begin{trivlist}
\item[\hskip \labelsep {\bfseries #1}]}{\end{trivlist}}

\newenvironment{rlist}{\begin{list}{$({\rm \roman{enumi}})$}
{\usecounter{enumi} \setlength{\rightmargin}{10pt}
\setlength{\leftmargin}{40pt} \setlength{\itemsep}{2pt}
\setlength{\parsep}{0pt} \setlength{\labelwidth}{40pt}}}{\end{list}}

\def\eq#1{{\rm(\ref{#1})}}

\def\N{{\mathbb N}}

\def\R{{\mathbb R}}

\def\GL{\mathbin{\rm GL}}

\def\SU{\mathop{\rm SU}}

\def\G2{\mathop\textrm{G}_2}

\def\w{\wedge}

\def\Nhat{\hat{N}}

\begin{document}

\title{Desingularization of Coassociative 4-folds with Conical Singularities}
\author{\textsc{Jason Dean Lotay}\\ University College\\ Oxford}
\date{}

\maketitle

\begin{center}
{\large\bf Abstract}
\end{center}

\noindent Suppose that $N$ is a compact \emph{coassociative} 4-fold
 with a \emph{conical singularity} in a 7-manifold $M$, with a $\G2$ structure given by a closed 3-form.  
We construct a smooth family, 
$\{N^\prime(t)\,:\,t\in(0,\tau)\}$ for some $\tau>0$, of 
compact, nonsingular, coassociative 4-folds in $M$ which converge to $N$ in the sense of currents, in geometric measure theory,
as $t\rightarrow 0$.  This realisation of \emph{desingularizations} of $N$ is achieved by gluing in
an \emph{asymptotically conical} coassociative 4-fold in $\R^7$, dilated by $t$, then deforming the resulting compact 4-dimensional 
submanifold of $M$ to the required coassociative 4-fold.  

\section{Introduction}

Harvey and Lawson, as part of their study of calibrated geometries \cite{HarLaw}, defined certain calibrated 4-dimensional
submanifolds of $\R^7$ known as \emph{coassociative} 4-folds and gave a family of $\SU(2)$-invariant examples.  Coassociative 
submanifolds are intimately linked to the imaginary octonions and thus to their automorphism group, $\G2$.  Details of the definition 
of these submanifolds and the means to extend the definition to more general 7-manifolds can be found in $\S$\ref{s2}.

Coassociative 4-folds with conical singularities, or CS coassociative 4-folds, were introduced and studied in the author's paper \cite{Lotay2},
which paid particular attention to their deformation theory. 
Asymptotically conical (AC) coassociative 4-folds in $\R^7$, and their deformations, were studied in the author's article
\cite{Lotay1}.  The results and techniques developed during these studies have provided the impetus for the research detailed here.

\subsection{Motivation}

The primary source of inspiration for the ideas in this paper, other than the author's earlier research, is the work of Joyce \cite{Joyce1}, \cite{Joyce2}, \cite{Joyce3},
\cite{Joyce4} and \cite{Joyce5} on special Lagrangian submanifolds with conical singularities, in particular \cite{Joyce3} on
desingularization.  There are a number of analogies between Joyce's work and the material in this article: in particular,
SL and coassociative submanifolds are both defined by the vanishing of a closed differential form and nearby deformations of these
submanifolds can be identified with graphs of
small forms.  This means that many of the analytic techniques and results, though not the same, 
are similar in style and content.  However, there are some considerable
discrepancies: notably, 
 Joyce is able to focus his study on functions whereas we must stick to forms, adding complexity to our problem.  

Another piece of work which motivates the study in this paper is by Kovalev \cite{Kovalev}, in which he constructs 
fibrations for a family of compact 7-manifolds with holonomy $\G2$ 
using coassociative $K3$ surfaces, the singular fibres of which are coassociative 4-folds with conical
singularities.  

\subsection{Summary}

We begin in $\S$\ref{s2} with the definition of coassociative 4-folds and $\varphi$-closed 7-manifolds, 
which are 7-manifolds with a $\G2$ structure given by a closed 3-form.  We continue by describing 
CS and AC coassociative 4-folds in $\S$\ref{s3}.  This section also contains material on weighted Banach spaces and elliptic operators
acting between such spaces, which provides the analytic framework for our study.

In $\S$\ref{s4} we construct, for all $t\in(0,\tau)$, a suitable compact (nonsingular) 4-dimensional
submanifold $\tilde{N}(t)$ in a $\varphi$-closed 7-manifold $M$, 
from a CS coassociative 4-fold $N$ in $M$ and an AC coassociative 4-fold $A$ in $\R^7$, 
such that $\tilde{N}(t)\rightarrow N$ as $t\rightarrow 0$ in the sense of currents.
Having presented this construction, we describe a carefully chosen \emph{tubular neighbourhood} of $\tilde{N}(t)$ in $\S$\ref{s5}, so that
we may characterise nearby deformations of $\tilde{N}(t)$ as the graphs of small forms on $\tilde{N}(t)$.
We then, in $\S$\ref{s6}, define a \emph{deformation map} $F_t$ on $\tilde{N}(t)$, whose kernel corresponds locally to coassociative deformations of $\tilde{N}(t)$.

In $\S$\ref{sei} we shift our attention to proving important analytic results on $\tilde{N}(t)$;
namely \emph{Sobolev embedding inequalities}
and \emph{regularity} results. 
Applying these results and further estimates, 
we define, in $\S$\ref{s7}, a \emph{contraction map} on a particular Banach space of forms on
$\tilde{N}(t)$.  The form fixed by this contraction is in the kernel of $F_t$.  In $\S$\ref{s9}, we
prove that this form is in fact smooth and thus corresponds to a coassociative deformation of $\tilde{N}(t)$.  

We conclude
with our key result, Theorem \ref{s9thm1}, from which we may deduce our main theorem, presented below.

\begin{thm}\label{s1thm1}  
Let $N$ be a coassociative 4-fold with $s$ distinct conical singularities $z_1,\ldots,z_s$ in a $\varphi$-closed 7-manifold $M$, 
as given in Definition \ref{singdfn2},
where the cone at $z_i$ is $C_i$ for $i=1,\ldots,s$.  Suppose that the link $\Sigma_i$ of the cone $C_i$ satisfies $b^1(\Sigma_i)=0$
for $i=1,\ldots,s$.
  Let $A_1,\ldots,A_s$ be coassociative 
4-folds in $\R^7$ such that $A_i$ is asymptotically conical to $C_i$ with rate $\lambda_i<-2$ for $i=1,\ldots,s$, as defined in Definition 
\ref{dt1dfn3}.  There exist $\tau>0$ 
and a smooth family, $\{N^{\prime}(t):t\in(0,\tau)\}$, of compact coassociative 4-folds in $M$ such that $N^{\prime}(t)\rightarrow N$ as
$t\rightarrow 0$, in the sense of currents in geometric measure theory.  Moreover, for all $t\in(0,\tau)$, 
$N^{\prime}(t)$ is diffeomorphic to the union of a compact subset of $N$, $C_i$ and a compact
subset of $A_i$ for $i=1,\ldots,s$.
\end{thm}

\begin{remark}
The conditions $b^1(\Sigma_i)=0$ and $\lambda_i<-2$ mean that this result is most closely analogous to \cite[Theorem 7.10]{Joyce3};
 that is, the main SL desingularization result in the \emph{unobstructed} case.
\end{remark}

\noindent To help the reader understand the nature of the convergence in this theorem, we briefly present 
some ideas from geometric measure theory.  

\begin{dfn}\label{s1dfn1}
Let $P$ be either a compact, oriented, 4-dimensional submanifold, or a 4-fold with conical singularities as defined in 
Definition \ref{singdfn2},
 in a 
7-manifold $M$.   Let $\mathcal{V}$ be the
space of smooth compactly supported 4-forms on $M$.  In the language of 
geometric measure theory, $P$ is an example of a \emph{oriented rectifiable set}.  Therefore, one may define, for 
$\xi\in\mathcal{V}$, the operator
$$\xi\mapsto \int_P\xi.$$
In this way, $P$ may be viewed as a 4-dimensional \emph{current}; that is, an element of the dual space $\mathcal{V}^*$.
Abusing notation slightly, we say that a sequence $(P_j)$ in $\mathcal{V}^*$ \emph{converges} to $P\in\mathcal{V}^*$ if it converges in the \emph{weak topology} on  $\mathcal{V}^*$;
i.e. 
$$P_j\rightarrow P\quad\text{if and only if}\quad \int_{P_j}\xi\rightarrow\int_P\xi\;\,\text{for all $\xi\in\mathcal{V}$.}$$
A good reference for this material is \cite[p. 30-34 \& p. 42-43]{Morgan}.
\end{dfn}

\begin{notes}\begin{itemize}\item[]
\item[(a)] Manifolds are taken to be nonsingular and submanifolds to
be embedded, for convenience, unless stated otherwise.
\item[(b)] We use the convention that the natural numbers
$\N=\{0,1,2,\ldots\}$.\end{itemize} \end{notes}

\section{Coassociative 4-folds}\label{s2}

\subsection{Basic definitions}

The key to defining coassociative 4-folds lies with the introduction
of a distinguished 3-form on $\R^7$.

\begin{dfn}\label{s2dfn1} Let $(x_1,\ldots,x_7)$ be coordinates on $\R^7$ and write
$d{\bf x}_{ij\ldots k}$ for the form $dx_i\w dx_j\w\ldots\w dx_k$.
Define a 3-form $\varphi_0$ by:
\begin{equation}\label{s2eq1}
\varphi_0 = d{\bf x}_{123}+d{\bf x}_{145}+d{\bf x}_{167}+d{\bf
x}_{246}- d{\bf x}_{257}-d{\bf x}_{347}-d{\bf x}_{356}.
\end{equation}
The 4-form $\ast\varphi_0$, where $\varphi_0$ and $\ast\varphi_0$
are related by the Hodge star, is given by:
\begin{equation*}\label{s2eq2}
\ast\varphi_0 = d{\bf x}_{4567}+d{\bf x}_{2367}+d{\bf
x}_{2345}+d{\bf x}_{1357}-d{\bf x}_{1346}-d{\bf x}_{1256}-d{\bf
x}_{1247}.
\end{equation*}
\end{dfn}

\vspace{-20pt}

\noindent Our choice of expression \eq{s2eq1} for $\varphi_0$
follows that of \cite[Chapter 10]{Joy1}.  This form is sometimes
known as the $\G2$ 3-form because the Lie group $\G2$ is the
subgroup of $\GL(7,\R)$ preserving $\varphi_0$.

\begin{dfn}\label{s2dfn2} A 4-dimensional submanifold $P$ of\/ $\R^7$
is coassociative if and only if $\varphi_0|_P\equiv 0$ and
$*\varphi_0|_P>0$.
\end{dfn}

\noindent This definition is not standard but is equivalent to the
usual definition in the language of \emph{calibrated geometry} by
\cite[Proposition IV.4.5 \& Theorem IV.4.6]{HarLaw}.

\begin{remark}
The condition $\varphi_0|_P\equiv 0$ forces $*\varphi_0$ to be a
non-vanishing 4-form on $P$.  Thus, the positivity of $*\varphi_0|_P$
is equivalent to a choice of orientation on $P$.
\end{remark}

So that we may describe coassociative submanifolds of more general
7-manifolds, we make two definitions following \cite[p. 7]{Bryant4}
and \cite[p. 243]{Joy1}.

\begin{dfn}\label{s2dfn3}
Let $M$ be an oriented 7-manifold.  For each $x\in M$ there exists
an orientation preserving isomorphism $\iota_x:T_xM\rightarrow\R^7$.
Since $\text{dim}\,\G2=14$, $\text{dim}\,\GL_+(T_xM)=49$ and
$\text{dim}\,\Lambda^3T^*_xM=35$, the $\GL_+(T_xM)$ orbit of
$\iota_x^*(\varphi_0)$ in $\Lambda^3T^*_xM$, denoted
$\Lambda^3_+T^*_xM$, is open.  A 3-form $\varphi$ on $M$ is
\emph{definite} 
if
$\varphi|_{T_xM}\in\Lambda^3_+T^*_xM$ for all $x\in M$.  Denote the
bundle of definite 3-forms $\Lambda^3_+T^*M$.  
\end{dfn}

\noindent Essentially, a definite 3-form is identified with the
$\G2$ 3-form on $\R^7$ at each point in $M$.  Therefore, to each
definite 3-form $\varphi$ we can uniquely associate a 4-form
$*\varphi$ and a metric $g$ on $M$ such that the triple
$(\varphi,*\varphi,g)$ corresponds to $(\varphi_0,*\varphi_0,g_0)$
at each point.  This leads us to our next definition.

\begin{dfn}\label{s2dfn4}
Let $M$ be an oriented 7-manifold, let $\varphi$ be a definite 3-form on $M$
 and let $g$ be the metric associated to $\varphi$.  We
call $(\varphi,g)$ a $\G2$ \emph{structure} on $M$. If $\varphi$ is
closed (or coclosed) then $(\varphi,g)$ is a \emph{closed} (or
\emph{coclosed}) $\G2$ structure.  A closed and coclosed $\G2$
structure is called \emph{torsion-free}.
\end{dfn}

\noindent Our choice of notation here agrees with \cite{Bryant4}.

\begin{remark}
By \cite[Lemma 11.5]{Salamon}, $(\varphi,g)$ is a torsion-free $\G2$ structure on $M$ if and only if 
the holonomy of $g$ is contained in $\G2$.
\end{remark}

\begin{dfn}\label{s2dfn5} Let $M$ be an oriented
7-manifold endowed with a $\G2$ structure $(\varphi,g)$, denoted
$(M,\varphi,g)$. We say that $(M,\varphi,g)$ is a
\emph{$\varphi$-closed}, or \emph{$\varphi$-coclosed}, 7-manifold if
$(\varphi,g)$ is a closed, respectively coclosed, $\G2$ structure.
If $(\varphi,g)$ is torsion-free, we call $(M,\varphi,g)$ a $\G2$
\emph{manifold}.
\end{dfn}

We are now able to complete our definitions regarding coassociative
4-folds.

\begin{dfn}\label{s2dfn6}
A 4-dimensional submanifold $P$ of\/ $(M,\varphi,g)$ is
coassociative if and only if\/ $\varphi|_P\equiv 0$ and
$*\varphi|_P>0$.
\end{dfn}

\begin{remark} Though we may define coassociative 4-folds with respect to any $\G2$ structure, for 
deformation theory and related results to hold we need it to be closed, but not necessarily coclosed.
Therefore, we shall work with $\varphi$-closed 7-manifolds for greatest useful generality.
\end{remark}

\subsection{Elementary results for 4-folds}

Since we will be concerned with self-dual 2-forms on 4-dimensional
manifolds with respect to different metrics, we make the following
definition.

\begin{dfn}\label{s2dfn7}
Let $(P,g_P)$ be an oriented, 4-dimensional, Riemannian manifold.  Denote the self-dual and
anti-self-dual 2-forms on $P$, with respect to $g_P$,
by $(\Lambda^2_+)_{g_P}T^*P$ and $(\Lambda^2_-)_{g_P}T^*P$.  Recall that if $P$ is compact,
there is a natural indefinite real-valued form on $\Lambda^2T^*P$, essentially given by integrating
 wedge product, which is positive definite on $(\Lambda^2_+)_{g_P}T^*P$.  We shall refer to this form as the \emph{intersection form}.
\end{dfn}

The next result, \cite[cf.$\!$ Proposition 4.2]{McLean},
is invaluable in describing the deformation theory of coassociative
4-folds.

\begin{prop}\label{s2prop2}
Let $P$ be a coassociative 4-fold in a 7-manifold $M$ with $\text{\emph{G}}_2$ structure 
$(\varphi,g)$. There is an
isometric isomorphism between the normal bundle $\nu(P)$ of $P$ in $M$ and
$(\Lambda^2_+)_{g|_P}T^*P$ given by $\jmath_P:v\mapsto
(v\cdot\varphi)|_{TP}$.
\end{prop}

We can generalise Proposition \ref{s2prop2} to 4-dimensional submanifolds of $(M,\varphi,g)$ which
are not necessarily coassociative.

\begin{prop}\label{s2prop4}
Let $P$ be a compact, oriented, 4-dimensional submanifold of a 7-manifold 
$M$ with $\text{\emph{G}}_2$ structure $(\varphi,g)$.
  Define $\jmath_P:\nu(P)\rightarrow\Lambda^2T^*P$ by
$\jmath_P:v\mapsto(v\cdot\varphi)|_{TP}$.  There exists a constant $\epsilon_{\varphi}>0$, independent of $P$, 
such that if $\|\varphi|_P\|_{C^0}<\epsilon_{\varphi}$ then there exists an orthogonal splitting
$$\Lambda^2T^*P=\Lambda^2_{\jmath+}T^*P\oplus\Lambda^2_{\jmath-}T^*P,$$
where $\nu(P)\stackrel{\jmath_P}{\cong}\Lambda^2_{\jmath+}T^*P$, for which the intersection form is positive definite on $\Lambda^2_{\jmath+}T^*P$.  
There also exists a unique metric
$g_P$ on $P$ 
such that
$$\Lambda^2_{\jmath+}T^*P=(\Lambda^2_+)_{g_P}T^*P,\quad 
\Lambda^2_{\jmath-}T^*P=(\Lambda^2_-)_{g_P}T^*P\;\,\text{and}\;\,*\varphi|_P=\text{\emph{vol}}_{g_P} .$$
\end{prop}

\begin{proof}
It is clear that $\jmath_P$ is a well-defined map from $\nu(P)$ into
$\Lambda^2T^*P$.  Further, from analysing the proof of
\cite[Proposition 4.2]{McLean} we deduce that $\jmath_P$ is
injective, irrespective of whether $\varphi$ vanishes on $P$ or not.
Defining $\Lambda^2_{\jmath+}T^*P=\jmath_P(\nu(P))$, we have a
3-dimensional subbundle of $\Lambda^2T^*P$ which is isomorphic to $\nu(P)$
via $\jmath_P$.  Moreover, the intersection form is positive definite on $\Lambda^2_{\jmath+}T^*P$, since this is an open condition and
holds when $\varphi|_P\equiv0$. Thus, we may construct an orthogonal splitting as claimed by defining
$$\Lambda^2_{\jmath_-}T^*P=\{\alpha\in\Lambda^2T^*P\,:\,\alpha\w\beta=0\;\,\text{for all $\beta\in\Lambda^2_{\jmath_+}T^*P$}\}.$$  

The existence of $g_P$ is a direct consequence of the
discussion in \cite[$\S$1.1.5]{DonKron}: in summary, for any 3-dimensional subbundle of $\Lambda^2T^*P$ which is positive with respect to
the intersection form, there exists a unique conformal structure on $P$ for which it is the self-dual subbundle.  Therefore the orthogonal
splitting of $\Lambda^2T^*P$ defines a unique conformal class of metrics on $P$ so that $\Lambda^2_{\jmath+}T^*P$ is the bundle of
self-dual 2-forms.  We then fix $g_P$ uniquely within this conformal
class by requiring the condition stated on its volume form.
\end{proof}

We finally prove two useful elementary results.  We denote Sobolev spaces by $L^p_k$ for $p\geq 1$ and $k\in\N$, 
as in \cite[$\S$1.2]{Joy1}.

\begin{prop}\label{s2prop3}
If $(P,g_P)$ is an oriented, compact, 4-dimensional, Riemannian manifold, $p\geq 1$ and $k\in\N$,
$$W^p_k=\{\gamma\in L^p_{k}(\Lambda^3T^*P)\,:\,\text{$\gamma$ is exact}\}=d\Big(L^p_{k+1}\big((\Lambda^2_+)_{g_P}T^*P\big)\Big).$$
\end{prop}

\begin{proof}  
Let $\gamma\in W^p_k$.
By Hodge theory, there exists a 2-form $\alpha$ in $L^p_{k+1}$ such 
that $d\alpha=\gamma$.  Using Hodge theory again, we may write $\alpha$ as the sum of an exact, a coexact and harmonic form, all in $L^p_{k+1}$.  
Since exact and harmonic forms are closed, there exists a 1-form $\xi$ and 
a closed 2-form $\eta$ such that $\alpha=*d\xi+\eta$.  Noting that $d(*d\xi+d\xi)=d*d\xi=d\alpha$, we see that
$\beta=*d\xi+d\xi$ is a self-dual 2-form in $L^p_{k+1}$ such that $d\beta=d\alpha=\gamma$.
\end{proof}

\begin{cor}\label{s2cor1} Use the notation of Proposition \ref{s2prop3}.
Let 
\begin{align*}
V^p_{k+1}=\big\{\alpha\in L^p_{k+1}\big((\Lambda^2_+)_{g_P}T^*P\big)\,&:\,\langle\alpha,\beta\rangle_{L^2}=0\\
&\,\;\;\,\text{for all closed}\,\,\beta\in C^{\infty}\big((\Lambda^2_+)_{g_P}T^*P\big)\big\}.\end{align*}
The map $d:V^p_{k+1}\rightarrow W^p_k$ is a linear isomorphism of Banach spaces and 
therefore has a linear inverse 
$d^{-1}:W^p_k\rightarrow V^p_{k+1}$. 
\end{cor}

\begin{proof}
By Proposition \ref{s2prop3} and the choice of $V^p_{k+1}$, $d:V^p_{k+1}\rightarrow W^p_{k}$ is a linear bijection.  
Since $V^p_{k+1}$ and $W^p_{k}$ are closed subspaces of Banach spaces, 
the first because the $L^2$-orthogonality condition is a closed one and
the second because it is the image of a Fredholm map, they are themselves Banach spaces.
\end{proof}

\section{Conical singularities and asymptotically\\ conical behaviour}\label{s3}

We present here the background material necessary for our study,
which is lifted from the author's earlier papers \cite[$\S$3]{Lotay2} and \cite[$\S$2]{Lotay1}
with only cosmetic changes.

\subsection{CS coassociative 4-folds}\label{s3subs1}

Let $B(0;\epsilon_M)$ denote the open ball about $0$ in $\R^7$ with radius
$\epsilon_M>0$. 
We define a preferred choice of local coordinates on a $\varphi$-closed
7-manifold near a finite set of points, which is an analogue of one given for
almost Calabi--Yau manifolds in \cite[Definition 3.6]{Joyce1}.

\begin{dfn}\label{ch8s1dfn1}\label{coords}
Let $(M,\varphi,g)$ be a $\varphi$-closed 7-manifold and let
$z_1,\ldots,z_s\in M$ be distinct points. There exist a constant $\epsilon_M\in(0,1)$,
an open set $V_i\ni z_i$ in $M$ with $V_i\cap V_j=\emptyset$ for
$j\neq i$ and a diffeomorphism
$\chi_i:B(0;\epsilon_M)\rightarrow V_i$ with $\chi_i(0)=z_i$,
for $i=1,\ldots,s$, such that $\zeta_i=d\chi_i|_0:\R^7\rightarrow
T_{z_i}M$ is an isomorphism identifying the standard $\text{G}_2$
structure $(\varphi_0,g_0)$ on $\R^7$ with the pair
$(\varphi|_{T_{z_i}M},g|_{T_{z_i}M})$. We call the set
$\{\chi_i:B(0;\epsilon_M)\rightarrow V_i:i=1,\ldots,s\}$ a
\emph{$\text{\emph{G}}_2$ coordinate system near
${z_1,\ldots,z_s}$}.

We say that two $\text{G}_2$ coordinate systems near
$z_1,\ldots,z_s$, with maps $\chi_i$ and $\tilde{\chi}_i$ for
$i=1,\ldots,s$ respectively, are \emph{equivalent} if
$d\tilde{\chi}_i|_0=d\chi_i|_0=\zeta_i$ for all $i$.
\end{dfn}

\begin{dfn}\label{ch8s1dfn2}\label{singdfn2}
Let $(M,\varphi,g)$ be a $\varphi$-closed 7-manifold, let $N\subseteq M$
be compact and connected and let $z_1,\ldots,z_s\in N$ be distinct. 
Let $\{\chi_i:B(0;\epsilon_M)\rightarrow
V_i:i=1,\ldots,s\}$ be a $\text{G}_2$ coordinate system near
$z_1,\ldots,z_s$, as in Definition \ref{coords}.
We say that
$N$ is a 4-fold in $M$ with \emph{conical singularities at
$z_1,\ldots,z_s$ with rate $\mu$}, denoted a \emph{CS 4-fold},
if $\hat{N}=N\setminus\{z_1,\ldots,z_s\}$ is a (nonsingular)
4-dimensional submanifold of $M$ and there exist constants
$0<\epsilon<\epsilon_M$ and $\mu\in(1,2)$, a compact 3-dimensional
Riemannian submanifold $(\Sigma_i,h_i)$ of
$\mathcal{S}^6\subseteq\R^7$, where $h_i$ is the restriction of the
round metric on $\mathcal{S}^6$ to $\Sigma_i$, an open set $U_i\ni
z_i$ in $N$ with $U_i\subseteq V_i$ and a smooth map
$\Phi_i:(0,\epsilon)\times\Sigma_i\rightarrow
B(0;\epsilon_M)\subseteq\R^7$, for $i=1,\ldots,s$, such that
 $\Psi_{i}=\chi_i\circ\Phi_i:(0,\epsilon)\times\Sigma_i\rightarrow
 U_i\setminus\{z_i\}$ is a diffeomorphism, and
 $\Phi_i$ satisfies 
\begin{equation}\label{ch8s1eq1a}
\Phi_i(r_i,\sigma_i)-\iota_i(r_i,\sigma_i)\in \big(T_{r_i\sigma_i}\iota_i(C_i)\big)^{\perp}\qquad 
\text{for all $(r_i,\sigma_i)\in (0,\epsilon)\times\Sigma_i$}
\end{equation} and
\begin{equation}\label{ch8s1eq1}
\big|\nabla^j_i\big(\Phi_i(r_i,\sigma_i)-\iota_i(r_i,\sigma_i)\big)\big|
=O\big(r_i^{\mu-j}\big)\qquad \text{for $j\in\N$ as $r_i\rightarrow 0$,}
\end{equation}
where $\iota_i(r_i,\sigma_i)=r_i\sigma_i\in B(0;\epsilon_M)$, $\nabla_i$
is the Levi--Civita connection of the cone metric
$g_i=dr_i^2+r_i^2h_i$ on $C_i=(0,\infty)\times\Sigma_i$ coupled with
partial differentiation on $\R^7$, and
 $|.|$ is calculated with respect to $g_i$.

 We call $C_i$ the \emph{cone} at the singularity $z_i$ and
$\Sigma_i$ the \emph{link} of the cone $C_i$.  We may write $N$ as
the disjoint union
$N=K_N\sqcup\bigsqcup_{i=1}^sU_i,$
where $K_N$ is compact.

If $\Nhat$ is coassociative in $M$, we say that $N$ is a \emph{CS
coassociative 4-fold}.
\end{dfn}

\begin{remark} If $N$ is a CS 4-fold, $\hat{N}$ is \emph{noncompact}.
\end{remark}

Suppose $N$ is a CS 4-fold at $z_1,\ldots,z_s$ with rate $\mu$
in $(M,\varphi,g)$ and use the notation of Definition
\ref{ch8s1dfn2}. The induced metric on $\Nhat$, $g|_{\Nhat}$, makes
$\Nhat$ into a Riemannian manifold. Moreover, it is clear from
\eq{ch8s1eq1} that, as long as $\mu<2$, the maps $\Psi_i$ satisfy
 \begin{equation*}\label{ch6s1eq2}
\big|\nabla_i^j\big(\Psi_i^*(g|_{\Nhat})-g_i\big)\big|=O\big(r_i^{\mu-1-j}\big)\qquad
\text{for $j\in\N$ as $r_i\rightarrow 0$.}
\end{equation*}
Consequently, the condition $\mu>1$ guarantees that the induced
metric on $\hat{N}$ genuinely converges to the conical metric on
$C_i$.

\begin{note} As shown on \cite[p. 6]{Lotay2}, since $\mu\in(1,2)$, Definition \ref{ch8s1dfn2} is
independent of the choice of $\text{G}_2$ coordinate system near the
singularities, up to equivalence.\end{note}

\begin{dfn}\label{ch6s1dfn4}\label{radiusfn}
Let $N$ be a CS coassociative 4-fold in a $\varphi$-closed 7-manifold $(M,\varphi,g)$ and use the notation of Definition \ref{singdfn2}. 
A \emph{radius function} on $\hat{N}$ is
a smooth map $\rho_N:\hat{N}\rightarrow (0,1]$ such that there exist positive constants
$c_1<1$ and $c_2>1$ with $c_1r_i<\Psi_i^*(\rho_N)<c_2r_i$ on
$(0,\epsilon)\times\Sigma_i$ for $i=1,\ldots,s$.  

It is clear how we may construct such a function.
\end{dfn}

\subsection{AC manifolds and submanifolds}\label{s3subs2}

\begin{dfn}
\label{dt1dfn3} Let $A$ be a closed (nonsingular) 4-dimensional
submanifold of $\R^7$. Then $A$ is \emph{asymptotically conical}
(AC), or an \emph{AC 4-fold}, \emph{with rate $\lambda$} if there
exist constants $\lambda<1$ and $R>0$, a compact subset $K_A$ of $A$,
a compact 3-dimensional Riemannian submanifold $(\Sigma,h)$ of
$\mathcal{S}^6\subseteq\R^7$, where $h$ is the restriction of the
round metric on $\mathcal{S}^6$ to $\Sigma$, and a diffeomorphism
$\Phi_A:(R,\infty)\times\Sigma\rightarrow A\setminus K_A$ satisfying 
\begin{equation}\label{dt1eq3a}
\Phi_A(r,\sigma)-\iota(r,\sigma)\in \big(T_{r\sigma}\iota(C)\big)^{\perp}\qquad\text{for all 
$(r,\sigma)\in (R,\infty)\times\Sigma$}\end{equation}
and
\begin{equation}
\label{dt1eq3}
\big|\nabla^j\big(\Phi_A(r,\sigma)-\iota(r,\sigma)\big)\big|=O\big(r^{\lambda-j}\big)
\qquad\text{for $j\in\N$ as $r\rightarrow\infty$,}
\end{equation}
\noindent where $\iota(r,\sigma)=r\sigma$, $\nabla$ is the
Levi--Civita connection of the cone metric $g_\text{cone}=dr^2+r^2h$
on $C=(0,\infty)\times\Sigma$ coupled with partial differentiation
on $\R^7$, and
 $|.|$ is calculated with respect to $g_\text{cone}$.

We call $\Sigma$ the \emph{link} of the cone $C$ as usual and say
that $A$ is asymptotically conical to $C$ for clarity when
necessary.
\end{dfn}

\begin{remark}
If $\lambda<0$ in the definition above, $A$ is AC with rate
$\lambda$ to a unique cone.
\end{remark}

In a similar manner to the CS case, we see that if $A$ is AC to $C$
with rate $\lambda$ then we can consider $A$ as a Riemannian
manifold with metric $g$ such that
\begin{equation*}
\label{dt1eq4}
\big|\nabla^j\big(\Phi_A^*(g)-g_{\text{cone}}\big)\big|=O\big(r^{\lambda-1-j}\big)\qquad\text{for
$j\in\N$ as $r\rightarrow\infty$},
\end{equation*}
by \eq{dt1eq3}, using the notation of Definition \ref{dt1dfn3}.  Thus, the constraint $\lambda<1$ ensures
that the metric on $A$ converges to $g_{\text{cone}}$.  

We also have the
analogous definition of a radius function in the AC case.

\begin{dfn}\label{ch6s1dfn2} Let $A$ be an AC 4-fold in $\R^7$ and use the notation of
Definition \ref{dt1dfn3}. A \emph{radius function}
$\rho_A:A\rightarrow [1,\infty)$ on $A$ is a smooth map such that
there exist positive constants $c_1<1$ and $c_2>1$ with
$c_1r<\Psi^*(\rho_A)<c_2r$ on $(R,\infty)\times\Sigma$.

It is again straightforward to see how one may construct
such a function.
\end{dfn}

\subsection{Weighted Banach spaces}\label{s3subs3}

 We use the notation and definition of the usual `unweighted' Banach spaces of forms as  
in \cite[$\S$1.2]{Joy1}; 
that is, Sobolev and H\"older spaces are denoted by $L^p_k$ and $C^{k,\,a}$ respectively,
where $p\geq 1$, $k\in\N$ and $a\in(0,1)$.  
Recall that, by the Sobolev Embedding Theorem, $L^p_{k}$ embeds continuously in $L^q_l$ if $l\leq k$ and
 $l-\frac{4}{q}\leq k-\frac{4}{p}$,
and $L^p_k$ embeds continuously in $C^{l,\,a}$ if $k-\frac{4}{p}\geq l+a$.
We also introduce the notation $C^k_{\text{loc}}$ for the space of forms $\xi$ such that $f\xi$ lies in $C^k$ for every
smooth compactly supported function $f$, and similarly define spaces $L^p_{k,\,\text{loc}}$ and $C^{k,\,a}_{\text{loc}}$.

We now define \emph{weighted} Banach spaces of forms as in
\cite[$\S$1]{Bartnik}.

\begin{dfn}\label{ch6s2dfn1}
%
%
Let $(P,g)$ be either the nonsingular part of a CS 4-fold or an AC 4-fold, as defined in Definitions \ref{singdfn2}
and \ref{dt1dfn3}.  Let $p\geq 1$, $k\in\N$, $\nu\in\R$ and let $\rho$
be a radius function on $P$, as described in Definitions \ref{radiusfn} and \ref{ch6s1dfn2}. The \emph{weighted Sobolev space}
$L_{k,\,\nu}^p(\Lambda^mT^*P)$ of $m$-forms $\xi$ on $P$ is the
subspace of
$L^p_{k,\,\text{loc}}(\Lambda^mT^*P)$ such that the norm
\begin{equation*}\label{ch6s2eq2}
\|\xi\|_{L_{k,\,\nu}^p}=\sum_{j=0}^k\left(\int_{P}
|\rho^{j-\nu}\nabla^j\xi|^p\rho^{-4} \,dV_g\right)^\frac{1}{p}
\end{equation*} is finite. 
Then $L_{k,\,\nu}^p(\Lambda^mT^*P)$ is a
Banach space and $L_{k,\,\nu}^2(\Lambda^mT^*P)$ is a Hilbert space.  
\end{dfn}

\begin{note} 
$L^p(\Lambda^mT^*P)=L_{0,\,-\frac{4}{p}}^p(\Lambda^mT^*P)$. 
\end{note}


\begin{dfn}\label{ch6s2dfn2} Let $P$ be either the nonsingular part of a CS 4-fold or an AC 4-fold, as defined in
Definitions \ref{singdfn2} and \ref{dt1dfn3}.  Let $\rho$
be a radius function on $P$, as given in Definition \ref{radiusfn} or \ref{ch6s1dfn2}, let $\nu\in\R$ and let $k\in\N$. The
\emph{weighted $C^k$-space} $C_{\nu}^{k}(\Lambda^mT^*P)$ of
$m$-forms $\xi$ on $P$ is the subspace of
$C^k_{\text{loc}}(\Lambda^mT^*P)$ such that the norm
$$\|\xi\|_{C_{\nu}^{k}}=\sum_{j=0}^k
\sup_{P}|\rho^{j-\nu}\nabla^j\xi|$$is finite. 
Then $C_{\nu}^{k}(\Lambda^mT^*P)$
is a Banach space. 
\end{dfn}

%

We shall need the analogue of the Sobolev Embedding Theorem for
weighted spaces, which is adapted from \cite[Lemma
7.2]{LockhartMcOwen} and \cite[Theorem 1.2]{Bartnik}.  It is
dependent on whether $P$ is a CS or AC 4-fold.
\begin{thm}[\textbf{Weighted Sobolev Embedding Theorem}]
\label{ch6s2thm1} Let $P$ be\\ either the nonsingular part of a CS 4-fold or an AC 4-fold, as given in Definitions \ref{singdfn2} and \ref{dt1dfn3},
and recall Definition \ref{ch6s2dfn1}.
Let $p,q\geq 1$, 
$\nu,\eta\in\R$ and $k,l\in\N$.
If $k\geq l$, $k-\frac{4}{p}\geq l-\frac{4}{q}$
and either
\begin{rlist}
\item $p\leq q$ and $\nu\geq\eta$ if $P$ is CS (or $\nu\leq\eta$ if $P$ is AC)
 or \item $p>q$ and $\nu>\eta$ if $P$ is CS (or $\nu<\eta$ if $P$ is AC),
\end{rlist}
there is a continuous embedding
$L_{k,\,\nu}^p(\Lambda^mT^*P)\hookrightarrow
L_{l,\,\eta}^q(\Lambda^mT^*P)$. 
\end{thm}


\subsection{Elliptic operators}\label{s3sub4}

We will mainly be concerned with elliptic operators acting on compact manifolds, but we shall have 
occasion to use elliptic regularity for operators on CS and AC 4-folds.  Let $P$ be the
nonsingular part of a CS 4-fold or an AC 4-fold, let $g$ be the metric on $P$, 
$p\geq 1$, $k\in\N$ and $\nu\in\R$.  We are interested in the map
\begin{equation}\label{d+d*}
d+d^*:L^p_{k+1,\,\nu}\big((\Lambda^2_+)_gT^*P\oplus\Lambda^4T^*P\big)\rightarrow L^p_{k,\,\nu-1}(\Lambda^3T^*P),
\end{equation}
where we identify the image of $d$ and $d^*$ from $(\Lambda^2_+)_gT^*P$.  This map is clearly elliptic and we can use 
the theory of \cite{LockhartMcOwen} and \cite{Lockhart}.  We start first with a result which follows from 
\cite[Theorems 1.1 and 6.1]{LockhartMcOwen}.

\begin{prop}\label{s3sub4prop1}
There exist countable discrete sets $\mathcal{D}_{\text{\emph{CS}}}$ and $\mathcal{D}_{\text{\emph{AC}}}$ of
real numbers such that \eq{d+d*} is Fredholm whenever either $\nu\notin\mathcal{D}_{\text{\emph{CS}}}$ or 
$\nu\notin\mathcal{D}_{\text{\emph{AC}}}$, depending on whether $P$ is either CS or AC.  
Moreover, $\mathcal{D}_{\text{\emph{CS}}}$ and $\mathcal{D}_{\text{\emph{AC}}}$ are determined by equations on 
$\Sigma_i$ and $\Sigma$ respectively, in the notation of Definitions \ref{singdfn2} and \ref{dt1dfn3}.
\end{prop}

From the work in \cite{Mazya}, which is in a more general setting, we deduce the following elliptic regularity result.

\begin{prop}\label{s3sub4prop3}
Let $P$ be the nonsingular part of a CS 4-fold or an AC 4-fold, as given in Definitions \ref{singdfn2} and \ref{ch6s1dfn2}, and let
$g$ be the metric on $P$.  Let $(\alpha,\beta)\in L^p_{1,\,\nu}\big((\Lambda^2_+)_gT^*P\oplus\Lambda^4T^*P\big)$ 
for some $p>1$ and $\nu\in\R$ and suppose that $d\alpha+d^*\beta\in L^p_{k,\,\nu-1}(\Lambda^3T^*P)$ for some $k\in\N$.  
Then $(\alpha,\beta)\in L^p_{k+1,\,\nu}\big((\Lambda^2_+)_gT^*P\oplus\Lambda^4T^*P\big)$ and there exists $c>0$, 
independent of $\alpha$ and $\beta$, such that
$$\|(\alpha,\beta)\|_{L^p_{k+1,\,\nu}}\leq c\left(\|d\alpha+d^*\beta\|_{L^p_{k,\,\nu-1}}+\|(\alpha,\beta)\|_{L^p_{1,\,\nu}}\right).$$
Moreover, if $\nu\notin\mathcal{D}_{\text{\emph{CS}}}$ or $\nu\notin\mathcal{D}_{\text{\emph{AC}}}$ as appropriate, 
where $\mathcal{D}_{\text{\emph{CS}}}$ and $\mathcal{D}_{\text{\emph{AC}}}$ are given by Proposition \ref{s3sub4prop1}, and 
$(\alpha,\beta)$ is $L^2$-orthogonal to the kernel of \eq{d+d*}, there exists $c^{\prime}>0$, again independent of
$\alpha$ and $\beta$, such that
\begin{equation*}\label{reg}
\|(\alpha,\beta)\|_{L^p_{k+1,\,\nu}}\leq c^{\prime}\|d\alpha+d^*\beta\|_{L^p_{k,\,\nu-1}}.
\end{equation*}
\end{prop}




We give a quick lemma now which we shall require later.  

\begin{lemma}\label{s3sub4lem1}
Using the notation of Definitions \ref{singdfn2} and \ref{dt1dfn3} and Proposition \ref{s3sub4prop1},
 $-2\notin\mathcal{D}_{\text{\emph{CS}}}$, or $-2\notin\mathcal{D}_{\text{\emph{AC}}}$, if and only if 
$b^1(\Sigma_i)=0$ for $i=1,\ldots,s$, or $b^1(\Sigma)=0$ respectively. 
\end{lemma}

\begin{proof}  For the AC case, this follows from the calculation preceding \cite[Proposition 5.2]{Lotay1},
which shows that $-2\in\mathcal{D}_{\text{AC}}$ if and only if there is a nonzero $(\alpha,\beta)\in C^{\infty}(\Lambda^2T^*\Sigma\oplus\Lambda^3T^*\Sigma)$ such that
$$d\alpha=-2\beta\quad\text{and}\quad d*\alpha+d^*\beta=0.$$
We see that $dd^*\beta=0$, so $\beta$ is harmonic and exact.  Hence $\beta=0$ by Hodge theory.  Therefore
$-2\in\mathcal{D}_{\text{AC}}$ if and only if there exists a nonzero closed and coclosed 2-form on $\Sigma$, i.e.
$b^2(\Sigma)=b^1(\Sigma)\neq 0$.  A similar calculation holds for the CS case.
\end{proof}

\section{Desingularization: stage 1}\label{s4}

\begin{note}
The notation introduced in this section shall be used throughout the remainder of the paper, though we
shall endeavour to remind the reader of key definitions when the need arises.
\end{note}

\subsection{The method and initial set-up}\label{s4subs1}

The first stage of desingularization is as follows.  We take a coassociative 4-fold in a $\varphi$-closed 7-manifold
$M$ with a single
conical singularity, modelled on a cone $C$, for convenience.  We then remove a small open neighbourhood
of the singularity and glue in a piece of a coassociative 4-fold in $\R^7$ which is asymptotically conical to $C$
and dilated by a factor $t$.  This process produces a compact (nonsingular), but
\emph{not necessarily coassociative}, 4-dimensional submanifold of $M$, which depends on $t$.  

\medskip

We thus need three sets of ingredients.

\begin{itemize}

\item[\textbf{(a)}] Let $(\Sigma,h)$ be a 3-dimensional Riemannian submanifold of $\mathcal{S}^6$, where $h$ is the restriction of the round metric
on $\mathcal{S}^6$ to $\Sigma$.  Furthermore, suppose $b^1(\Sigma)=0$.  Define $C=(0,\infty)\times\Sigma$, let $g_{\text{cone}}=dr^2+r^2h$ and
 let $\iota:C\rightarrow\R^7$ be the inclusion map,
$\iota(r,\sigma)=r\sigma$.  

\item[\textbf{(b)}] Let $N$ be a CS coassociative 4-fold in a $\varphi$-closed 7-manifold $(M,\varphi,g)$ with one
singular point $z$, with rate $\mu\in(1,2)$ and cone $C$, as defined in Definition \ref{singdfn2}.  
We thus require a $\G2$ coordinate system, as given in
Definition \ref{coords}, consisting of a single diffeomorphism $\chi:B(0;\epsilon_M)\rightarrow V$, where $\epsilon_M\in(0,1)$ is a constant
and $V$ is an open set in $M$ containing $z$.  Recall that $\hat{N}=N\setminus\{z\}$ and we have a constant $\epsilon\in(0,\epsilon_M)\subseteq(0,1)$ and
a smooth map $\Phi_N:(0,\epsilon)\times\Sigma\rightarrow B(0;\epsilon_M)$ satisfying \eq{ch8s1eq1a} and \eq{ch8s1eq1}.  
Moreover, there
exists a compact subset $K_N$ of $N$ such that
$z\in N\setminus K_N\subseteq V$ and $\hat{N}\setminus K_N$ is diffeomorphic to $(0,\epsilon)\times\Sigma$ via
$\Psi_N=\chi\circ\Phi_N$.

\item[\textbf{(c)}] Let $A$ be a coassociative 4-fold in $\R^7$ which is AC with rate
$\lambda<-2$ to $C$, as defined in Definition \ref{dt1dfn3}.  Then there is a compact subset
$K_A$ of $A$ such that $A\setminus K_A$ is diffeomorphic to $(R,\infty)\times\Sigma$ via $\Phi_A$, 
where $R>1$ is a constant and $\Phi_A$ satisfies \eq{dt1eq3a} and \eq{dt1eq3}.  
Finally, we let $t\in(0,\tau)$, where $\tau>0$ is small enough that $\tau^{-1}\epsilon>R$, and let 
$Y_A(t)=K_A\sqcup\Phi_A\big((R,t^{-1}\epsilon)\times\Sigma\big)$. It is clear that $tA$ is coassociative
and AC with rate $\lambda$ to $C$.

\end{itemize}

For our method to work, we must restrict our choices of $t$ by decreasing $\tau$.  We shall have to do this
a finite number of times during the course of this article before the conclusion.  At each juncture, we assume $\tau$ is reduced to 
satisfy the conditions necessary for the result to hold for all $t\in(0,\tau)$.  
This assumption will be encompassed in the simple statement that $\tau$ is ``sufficiently
small'', rather than cluttering the exposition with a sequence of decreasing constants $\tau_1$, $\tau_2$, $\tau_3$ etc.  

Our first result
constraining the size of $\tau$ is given below.  The proof is elementary and hence omitted.

\begin{lemma}\label{TA}
If $\tau$ is sufficiently small,
the sets $tY_A(t)$ and $t\iota\big((R,t^{-1}\epsilon)\times\Sigma\big)$ are contained
in $B(0;\epsilon_M)$.
\end{lemma}

\subsection{Constructing {\boldmath $\tilde{N}(t)$}}\label{s4subs2}

Armed with the set-up of $\S$\ref{s4subs1}, we proceed in defining, for each $t\in(0,\tau)$, a
compact nonsingular 4-fold $\tilde{N}(t)$ using $N$ and $tA$.

\begin{dfn}\label{s4dfn2} Use the notation of $\S$\ref{s4subs1}(a)-(c).
Let $f_{\text{inc}}:\R\rightarrow[0,1]$ be a smooth increasing
function such that
$$f_{\text{inc}}(x)=\left\{\begin{array}{lll} 0 & & \text{for $x\leq 0$,}\\
1&&\text{for $x\geq1$}\end{array}\right.$$ and
$f_{\text{inc}}(x)\in(0,1)$ for $x\in(0,1)$.  Let
$\nu\in(0,1)$, with $\nu>
\frac{3}{\mu+2}
$,
and choose $\tau$ sufficiently small so that
$0<tR<\frac{1}{2}t^{\nu}<t^{\nu}<\epsilon$, which is possible since $\mu>1$ and $\nu<1$.  Define a smooth map
$\Phi_{\tilde{N}(t)}:(tR,\epsilon)\times\Sigma\rightarrow B(0;\epsilon_M)$
by
\begin{align*}\Phi_{\tilde{N}(t)}(r,\sigma)&=t\big(1-f_\text{inc}( 
2t^{-\nu}r-1)\big)
\Phi_A(t^{-1}r,\sigma) 
+f_{\text{inc}}\left( 
2t^{-\nu}r-1\right)\Phi_N(
r,\sigma).
\end{align*} By construction the image of
$\Phi_{\tilde{N}(t)}$ is contained in $B(0;\epsilon_M)$ and can be considered
as a deformation of $\tilde{C}(t)=\iota\big((tR,\epsilon)\times\Sigma\big)$. 
Therefore, define
$\Psi_{\tilde{N}(t)}=\chi\circ\Phi_{\tilde{N}(t)}:(tR,\epsilon)\times\Sigma\rightarrow
V\subseteq M$ and let
$$\tilde{N}(t)=\chi(tK_A)\sqcup\Psi_{\tilde{N}(t)}\big((tR,\epsilon)\times\Sigma\big)\sqcup K_N.$$
Define
\begin{align*}
\tilde{N}_l(t)&=\chi(tK_A)\sqcup\Psi_{\tilde{N}(t)}\big((tR,\textstyle\frac{1}{2}\displaystyle t^{\nu})\times\Sigma\big),\\
\tilde{N}_m(t)&=\Psi_{\tilde{N}(t)}\big([\textstyle\frac{1}{2}\displaystyle t^{\nu},t^{\nu}]\times\Sigma\big)\,\;\text{and}\\
\tilde{N}_u(t)&=\Psi_{\tilde{N}(t)}\big((t^{\nu},\epsilon)\times\Sigma\big)\sqcup
K_N.\end{align*} 
It is clear that $\tilde{N}(t)$ is a compact
(nonsingular) 4-dimensional submanifold of $M$, which is the disjoint union of $\tilde{N}_l(t)$, $\tilde{N}_m(t)$ and
$\tilde{N}_u(t)$.  Moreover, $\tilde{N}_l(t)\subseteq\chi\big(tY_A(t)\big)$ is an open subset of $\tilde{N}(t)$, 
$\tilde{N}_m(t)$ is a compact subset of $\tilde{N}(t)$ and $\tilde{N}_u(t)\subseteq\hat{N}$ is an open subset of $\tilde{N}(t)$.
\end{dfn}

\begin{remark}
The choice of $\nu$ shall become clearer later on, but the idea is that there is a balancing act between the 
`interpolation' region being small but not too tight, 
and the number $t^{\nu}$ being small whilst $t^{\nu-1}$ is large.
\end{remark}

\begin{note}
The definition of $\tilde{N}(t)$ is such that, as $t\rightarrow 0$, 
$\tilde{N}(t)\setminus\tilde{N}_u(t)$ collapses to the singularity $z$
whilst $\tilde{N}_u(t)$ expands to equal $\hat{N}$.  Thus, $\tilde{N}(t)\rightarrow N$ in the sense of currents as $t\rightarrow 0$.  
\end{note}

Though we defer the detailed calculations estimating the modulus of $\varphi|_{\tilde{N}(t)}$ until $\S$\ref{s7}, we require
the following weaker result now.

\begin{lemma}\label{smallphi}
Let $\epsilon_{\varphi}>0$ be the constant given in Proposition \ref{s2prop4} applied to $(M,\varphi,g)$.  
If $\tau$ is sufficiently small,
$\big\|\varphi|_{\tilde{N}(t)}\big\|_{C^0}<\epsilon_{\varphi}$.  
\end{lemma}

\begin{proof}
Later, in Proposition \ref{s7prop1}, we estimate $\varphi|_{\tilde{N}(t)}$.  
From the calculations in the 
proof of that proposition we deduce that 
$$\Big\|\chi^*(\varphi)|_{\chi^{-1}\big(\tilde{N}_l(t)\big)}\Big\|_{C^0}\rightarrow 0\quad\text{and}\quad
\Big\|\Psi_{\tilde{N}_t}^*\Big(\varphi|_{\tilde{N}_m(t)}\Big)\Big\|_{C^0}\rightarrow 0$$
as $t\rightarrow 0$, using the Euclidean metric on $\R^7$ for the first term and the conical metric on $[\frac{1}{2}t^{\nu},t^{\nu}]
\times\Sigma$ for the second.  Note that these calculations do \emph{not} require that $\tilde{N}(t)$ satisfies Proposition 
\ref{s2prop4}, although we do use that fact in the proof of Proposition \ref{s7prop1}.  
Clearly $\varphi|_{\tilde{N}_u(t)}\equiv 0$ as $\tilde{N}_u(t)\subseteq \hat{N}$.  Thus we quickly see
that $\|\varphi|_{\tilde{N}(t)}\|_{C^0}\rightarrow 0$ as $t\rightarrow 0$.
\end{proof}

In a similar vein to CS and AC 4-folds, we introduce a radius function on $\tilde{N}(t)$.  

\begin{dfn}\label{s4dfn3} Use the notation of $\S$\ref{s4subs1}(a)-(c) and Definition \ref{s4dfn2}.
A \emph{radius function} $\rho_t:\tilde{N}(t)\rightarrow [t,1]$ on $\tilde{N}(t)$ is a smooth map, with 
$\rho_t(x)=1$ for all $x\in K_N$ and $\rho_t(x)=t$ for all $x\in\chi(tK_A)$, such that there exist constants $c_1<1$, $c_2>1$
and $c_3>0$, independent of $t$, with
$$c_1r<\Psi_{\tilde{N}(t)}^*(\rho_t)<c_2r\quad\text{and}\quad |\Psi_{\tilde{N}(t)}^*(d\rho_t)|\leq c_3$$ 
on $(tR,\epsilon)\times\Sigma$.
\end{dfn}

\noindent The existence of such a function on $\tilde{N}(t)$ is clear.  Given a radius function on $\tilde{N}(t)$ we can define 
weighted Banach spaces on it.   

\begin{dfn}\label{s4dfn4}  Use the notation of Definition \ref{s4dfn2} and let $\rho_t$ be a radius function on $\tilde{N}(t)$
as in Definition \ref{s4dfn3}. 
By Lemma \ref{smallphi}, Proposition \ref{s2prop4} is applicable to $\tilde{N}(t)$ so, in the notation of
that proposition, let $\tilde{g}(t)=g_{\tilde{N}(t)}$ and let $\nabla_t$ be the Levi--Civita connection of $\tilde{g}(t)$.  

For $k\in\N$ and $\eta\in\R$, define $C^k_{\eta,\,t}\big(\Lambda^mT^*\tilde{N}(t)\big)$ to be the subspace of
 $C^k_{\text{loc}}\big(\Lambda^mT^*\tilde{N}(t)\big)$ such that the norm
$$\|\xi\|_{C^k_{\eta,\,t}}=\sum_{j=0}^k\sup_{\tilde{N}(t)}|\rho_t^{j-\eta}\nabla_t^j\xi|,$$
calculated using $\tilde{g}(t)$, is finite.  These are Banach space which are analogues of the weighted $C^k$ spaces
defined in Definition \ref{ch6s2dfn2}.

For $p\geq 1$, $k\in\N$ and $\eta\in\R$, let $L^p_{k,\,\eta,\,t}\big(\Lambda^mT^*\tilde{N}(t)\big)$ be the subspace of
 $L^p_{k,\,\text{loc}}\big(\Lambda^mT^*\tilde{N}(t)\big)$ such that the norm
$$\|\xi\|_{L^p_{k,\,\eta,\,t}}=\sum_{j=0}^k\left(\int_{\tilde{N}(t)}|\rho_t^{j-\eta}\nabla_t^j\xi|^p
\rho_t^{-4}
dV_{\tilde{g}(t)}\right)^{\frac{1}{p}}$$
is finite.  These are Banach spaces analogous to the weighted Sobolev spaces given in Definition \ref{ch6s2dfn1}.
\end{dfn}

\begin{notes}
\begin{itemize}\item[]
\item[(a)]
The spaces $C^k_{\eta,\,t}$ and 
$L^p_{k,\,\eta,\,t}$ are Lipschitz equivalent to the corresponding `unweighted' spaces, $C^k$ and $L^p_k$,
 because $\tilde{N}(t)$ is compact. 
However, the Lipschitz constants depend on $t$.  This is why we introduce these `weighted' spaces: to manage the $t$ dependence of the norm.
\item[(b)]
It is clear there is a continuous embedding of $L^p_{2,\,1,\,t}$ into $C^1_{1,\,t}$ for $p>4$ by analogy
with the Weighted Sobolev Embedding Theorem (Theorem \ref{ch6s2thm1}).  Moreover, by studying the $t$ dependence of the norms,
we see that the embedding constant is independent of $t$.  This fact may be proved using the same method as
Proposition \ref{s5prop8} so we omit it.  

The idea is to consider small balls on 
$\tilde{N}(t)$ on which we have the usual embedding inequality, with
 embedding constant multiplied by a factor determined by the radius of the ball.  The weighted norms take care of these
 radius factors.  Thus, using a finite open cover by balls 
 given by the compactness of $\tilde{N}(t)$ yields the global embedding inequality.  
\end{itemize}   
\end{notes}

\subsection{Cohomology of {\boldmath $\tilde{N}(t)$}}

We show how we can describe the cohomology of $\tilde{N}(t)$ using the constituent parts $\hat{N}$ and $A$.  
From standard results in algebraic topology, 
considering $\hat{N}$ as the interior of a manifold with boundary $\Sigma$, there exist homomorphisms 
$\iota_m^N$, $j_m^N$ and $\partial_m^N$ such that the following sequence is exact:
\begin{equation}\label{Nseq}
\cdots\longrightarrow H^m_{\text{cs}}(\hat{N})
\,{\buildrel\iota^N_m\over\longrightarrow}\,
H^m_{\text{dR}}(\hat{N})\,{\buildrel j^N_m\over\longrightarrow}\,
H^m_{\text{dR}}(\Sigma)\,{\buildrel
\partial^N_m\over\longrightarrow}\,
H_{\text{cs}}^{m+1}(\hat{N})\longrightarrow\cdots.
\end{equation}
Moreover, $\iota_m^N$ is the inclusion map.

\begin{dfn}\label{IN}
Using the notation in \eq{Nseq}, let $\mathcal{I}^N=\iota^N_2\big(H^2_{\text{cs}}(\hat{N})\big)$.  By
\cite[Example (0.16)]{Lockhart}, the space
$$\mathcal{H}^2(N)=\{\xi\in L^2(\Lambda^2T^*\hat{N})\,:\,d\xi=d^*\xi=0\}$$
is isomorphic to $\mathcal{I}^N$ via $\xi\mapsto [\xi]$.  Note that $\mathcal{H}^2(N)$ consists of smooth forms by
elliptic regularity.  

Define a cup product on $\mathcal{I}^N$ as follows.  Let $\alpha$ and $\beta$ be elements of $\mathcal{I}^N$ and let
$\xi$ and $\eta$ be closed compactly supported 2-forms on $\hat{N}$ such that $\iota^N_2([\xi])=\alpha$ and
$\iota^N_2([\eta])=\beta$.  Then define
$$\alpha\cup\beta=\int_{\hat{N}}\xi\w\eta.$$
In \cite[Definition 8.1]{Lotay2}, the author showed that this gives a well-defined symmetric topological product on
$\mathcal{I}^N$.  Furthermore,
$$\mathcal{H}^2_+(N)=\mathcal{H}^2(N)\cap C^{\infty}(\Lambda^2_+T^*\hat{N})$$
is isomorphic to the positive subspace of $\mathcal{I}^N$ with respect to the cup product.
\end{dfn}

Again using standard theory, there exist homomorphisms $\iota_m^A$, $j_m^A$ and $\partial_m^A$, with $\iota_m^A$ the obvious 
inclusion map, such that the following is exact:
\begin{equation}\label{Aseq}
\cdots\longrightarrow H^m_{\text{cs}}(A)
\,{\buildrel\iota^A_m\over\longrightarrow}\,
H^m_{\text{dR}}(A)\,{\buildrel j^A_m\over\longrightarrow}\,
H^m_{\text{dR}}(\Sigma)\,{\buildrel
\partial^A_m\over\longrightarrow}\,
H_{\text{cs}}^{m+1}(A)\longrightarrow\cdots.
\end{equation}

\begin{dfn}\label{IA}
Using the notation in \eq{Aseq}, let $\mathcal{I}^A=\iota^A_2\big(H^2_{\text{cs}}(A)\big)$.  By
\cite[Example (0.15)]{Lockhart}, the space
$$\mathcal{H}^2(A)=\{\xi\in L^2(\Lambda^2T^*A)\,:\,d\xi=d^*\xi=0\}$$
is isomorphic to $\mathcal{I}^A$ via $\xi\mapsto [\xi]$.  Note that $\mathcal{H}^2(A)$ consists of smooth forms by
elliptic regularity.    

We can again define a cup product on $\mathcal{I}^A$ and
$$\mathcal{H}^2_+(A)=\mathcal{H}^2(A)\cap C^{\infty}(\Lambda^2_+T^*A)$$
is isomorphic to the positive subspace of $\mathcal{I}^A$ with respect to the cup product.  

Note that this definition only requires $\lambda<1$.
\end{dfn}

Now we can consider $\tilde{N}(t)$ as the union of two open sets
 $$W_A=\chi(tK_A)\cup\Psi_{\tilde{N}(t)}\big((tR,t^{\nu-1})\times\Sigma\big)\quad\text{and}\quad 
W_N=K_N\cup\Psi_{\tilde{N}(t)}\big((\textstyle\frac{1}{2}t^{\nu-1},\epsilon)\times\Sigma\big).$$
Note that $W_A$ is diffeomorphic to $A$, $W_N$ is diffeomorphic to $\hat{N}$ 
and that $W_A\cap W_N$ is diffeomorphic to $C$.  Hence we can apply the Mayer--Vietoris Theorem to $\tilde{N}(t)$ 
and deduce that there exist homomorphism $\tilde{i}_m$, $\tilde{j}_m$ and $\tilde{\partial}_m$ such that the following
sequence is exact:
\begin{equation}\label{seq}
\cdots\!\longrightarrow H^m_{\text{dR}}\big(\tilde{N}(t)\big)
\,{\buildrel\tilde{i}_m\over\longrightarrow}\,
H^m_{\text{dR}}(A)\oplus H^m_{\text{dR}}(\hat{N})\,{\buildrel \tilde{j}_m\over\longrightarrow}\,
H^m_{\text{dR}}(\Sigma)\,{\buildrel
\tilde{\partial}_m\over\longrightarrow}\,
H_{\text{dR}}^{m+1}\big(\tilde{N}(t)\big)\longrightarrow\!\cdots.
\end{equation}

We now have a useful corollary.

\begin{cor}\label{seqcor} In the notation of $\S$\ref{s4subs1} and Definitions \ref{s4dfn2}, \ref{IN} and \ref{IA},
$$b^2_+\big(\tilde{N}(t)\big)= \text{\emph{dim}}\,\mathcal{H}^2_+(N)+\text{\emph{dim}}\,\mathcal{H}^2_+(A).$$
\end{cor}

\begin{proof}
Use the notation of equations \eq{Nseq}, \eq{Aseq} and \eq{seq}.  Since $b^1(\Sigma)=b^2(\Sigma)=0$ and the sequence 
\eq{seq} is exact, $H^2_{\text{dR}}(N)\cong H^2_{\text{dR}}(A)\oplus H^2_{\text{dR}}(\hat{N})$.  Moreover, 
as \eq{Nseq} and \eq{Aseq} are exact, $\iota_2^N$ and $\iota_2^A$ are isomorphisms, so
$\mathcal{I}^N\cong H^2_{\text{dR}}(\hat{N})$ and $\mathcal{I}^A\cong H^2_{\text{dR}}(A)$.  
The corollary follows from the observations in Definitions \ref{IN} and \ref{IA}.
\end{proof}

\section{Tubular neighbourhoods}\label{s5}

The main goal of this section, achieved in Proposition \ref{s5prop3}, is to construct a tubular neighbourhood of $\tilde{N}(t)$
which can be identified with a $C^1_{1,\,t}$-open set of small self-dual 2-forms on $\tilde{N}(t)$.  The key is to build this
neighbourhood from neighbourhoods of $C$, $N$ and $A$.

We start with an elementary result, pertinent to a general situation,
which is immediate from the proof of \cite[Chapter IV, Theorem
9]{Lang}.

\begin{thm}\label{s5thm1} Let $P$ be a closed embedded submanifold
of a Riemannian manifold $M^{\prime}$.  There exist an open subset $V_P$ of
the normal bundle $\nu(P)$ of $P$ in $M^{\prime}$, containing the zero
section, and an open set $S_P$ in $M^{\prime}$ containing $P$, such that the
exponential map $\exp|_{V_P}:V_P\rightarrow S_P$ is a
diffeomorphism.
\end{thm}

\begin{note} The proof of this result relies entirely on the observation
that $\exp|_{\nu(P)}$ is a local isomorphism upon the zero
section.\end{note}

\subsection{Corollaries of Theorem \ref{s5thm1}}\label{s5subs1}

Here we do the groundwork in constructing our neighbourhood of $\tilde{N}(t)$ by considering tubular neighbourhoods
of the cone $C$, the CS coassociative 4-fold $N$ and the AC coassociative 4-fold $A$.

To consider neighbourhoods of $C$ we first need a definition.

\begin{dfn} Recall the definition of $C$ in $\S$\ref{s4subs1}(a).
Let $v\in\nu(C)$, the normal bundle of $\iota(C)$ in $\R^7$. 
 For all $(r,\sigma)\in C$, $v|_{(r,\sigma)}\in\R^7$ and so we can consider the dilated vector
$s(v|_{(r,\sigma)})$ for $s>0$.  Thus, we can define an action of dilation on $\nu(C)$ by
$v\mapsto sv$ where $(sv)|_{(sr,\sigma)}=s(v|_{(r,\sigma)})$, $s>0$.
\end{dfn}

\begin{cor}\label{s5cor0} Recall the notation of $\S$\ref{s4subs1}(a). 
There exist a dilation-invariant open subset $V_C$ of the normal bundle $\nu(C)$ of $\iota(C)$ in $\R^7$,
containing the zero section, a dilation-invariant open set $S_C$ in $\R^7$ containing $C$, and a 
dilation-equivariant diffeomorphism $n_C:V_C\rightarrow S_C$ with $n_C(0)=C$.  Moreover, $V_C$ and $S_C$ grow with 
order $O(r)$ as $r\rightarrow 0$ and $r\rightarrow\infty$.
\end{cor}

\begin{proof}
Consider $P=\iota(\{1\}\times\Sigma)\subseteq\R^7$.  By Theorem \ref{s5thm1}, there exist
 an open subset $V_{\Sigma}$ of $\nu(P)$ and an open set $S_{\Sigma}$ in $\R^7$ which are
 diffeomorphic via the exponential map.  We can then define $S_C=\{rs\,:\,r>0, s\in S_{\Sigma}\}$,
 which contains $\iota(C)$, and similarly extend $V_{\Sigma}$ to $V_C$ in a dilation-invariant 
 way.  The diffeomorphism $n_C$ is then given naturally by $n_C(r\sigma,v)=r\exp(\sigma,v)$, where
 we identify $\nu_{r\sigma}(\iota(C))$ with $\nu_{\sigma}(P)$.  The final part is a direct consequence of
 the construction of the open sets described.
\end{proof}

Armed with this result, we can 
find tubular neighbourhoods related to our AC and CS coassociative 4-folds, which are conical around the 
cone-like parts of the coassociative 4-folds.

\begin{cor}\label{s5cor1} Use the notation of Definition \ref{ch6s2dfn2}, $\S$\ref{s4subs1}(a) and (c) and Corollary \ref{s5cor0}. 
Let $C_{\infty}=\iota\big((R,\infty)\times\Sigma\big)$. 
\begin{itemize}\item[\emph{(a)}] There exists a diffeomorphism $\upsilon_A:\nu(C_{\infty})\rightarrow
\nu(A\setminus K_A)$ which identifies the zero section in each bundle and maps $\nu_{r\sigma}(C_{\infty})$ 
to $\nu_{\Phi_A(r,\sigma)}(A\setminus K_A)$ for all $(r,\sigma)\in(R,\infty)\times\Sigma$.
\item[\emph{(b)}] There exist an open subset $V_A$ of the normal bundle
$\nu(A)$ of $A$ in $\R^7$, containing the zero section,
an open set $S_A$ in $\R^7$, containing $A$, and a diffeomorphism $n_A:V_A\rightarrow S_A$.
Moreover, $V_A$, $S_A$ and $n_A$ over $A\setminus K_A$ agree with $V_C$, $S_C$ and $n_C$ over $C_{\infty}$, where
we use $\upsilon_A$ in this identification, and $V_A$ may be chosen to be an open neighbourhood of the 
zero section in $C^1_1$.
\end{itemize}
\end{cor}

\begin{proof}
Recall from Corollary \ref{s5cor0} that $S_C$ is dilation-invariant and grows with order $O(r)$ as 
$r\rightarrow\infty$.  Therefore, by making $R$ and $K_A$ larger if necessary, we can ensure that 
$A\setminus K_A$ is contained in $S_C$, since $\Phi_A(r,\sigma)-r\sigma$
is order $O(r^{\lambda})$, where $\lambda<1$, by \eq{dt1eq3}.  
Using \eq{dt1eq3a}, since $(T_{r\sigma}C_{\infty})^{\perp}\cong\nu_{r\sigma}(C_{\infty})$, we have a diffeomorphism $\upsilon_A$ as claimed.  

Let $V_C^{\prime}=\{v\in\nu(C_{\infty})\,:\,v\in V_C\}$ and
define $V_A^{\prime}\subseteq\nu(A\setminus K_A)$ by $\upsilon_A(V_C^{\prime})$: this clearly contains the zero section over
$A\setminus K_A$.
We also define $S_A^{\prime}\subseteq S_C$
to be equal to $n_C(V_C^{\prime})$, which contains $A\setminus K_A$ and $C_{\infty}$ by the discussion above.
Thus $n_A^{\prime}=n_C\circ\upsilon_A^{-1}$ is a diffeomorphism between $V_A^{\prime}$
and $S_A^{\prime}$.
We then extend $V_A^\prime$,
$S_A^{\prime}$ and $n_A^{\prime}$ smoothly over the compact set $K_A$ to the sets and diffeomorphism we require.
Since $V_C$ grows with order $O(r)$ as $r\rightarrow\infty$, we can ensure $V_A$ is an open
neighbourhood in $C^1_1$.
\end{proof}

\begin{cor}\label{s5cor2}  Use the notation of Definition \ref{ch6s2dfn2}, $\S$\ref{s4subs1}(a) and (b) and Corollary \ref{s5cor0}.
Let
$C_0=\iota\big((0,\epsilon)\times\Sigma\big)$ and let $B=\Phi_N\big((0,\epsilon)\times\Sigma\big)$.
\begin{itemize}\item[\emph{(a)}] There exists a diffeomorphism $\upsilon_N:\nu(C_0)\rightarrow
\nu(B)$ which identifies the zero section in each bundle and maps $\nu_{r\sigma}(C_0)$ to $\nu_{\Phi_N(r,\sigma)}(B)$ 
for all $(r,\sigma)\in(0,\epsilon)\times\Sigma$.
\item[\emph{(b)}] There exist an open subset $V_B$ of the normal bundle
$\nu(B)$ of $B$ in $\R^7$, containing the zero section,
an open set $S_B$ in $B(0;\epsilon_M)$, containing $B$, and a diffeomorphism $n_B:V_B\rightarrow S_B$.
Moreover, $V_B$, $S_B$ and $n_B$ over $B$ agree with $V_C$, $S_C$ and $n_C$ over $C_0$, where
we use $\upsilon_B$ in this identification, and $V_B$ may be chosen to be an open neighbourhood of the 
zero section in $C^1_1$.
\end{itemize}
\end{cor}

\noindent The proof is almost identical to that of Corollary \ref{s5cor1} so we omit it.  We note we can ensure
that $V_B$ is an open neighbourhood of the zero section in $C^1_1$ since $V_C$ grows with order $O(r)$ as $r\rightarrow0$.

We can now combine the work above to give our next corollary to Theorem \ref{s5thm1},
which is our main building block in constructing a tubular neighbourhood of $\tilde{N}(t)$.

\begin{cor}\label{s5cor4}  Use the notation of $\S$\ref{s4subs1}(a)-(c), Definition \ref{s4dfn2} and Corollary \ref{s5cor0}.
Recall that $\tilde{C}(t)=\iota\big((tR,\epsilon)\times\Sigma\big)$ and let $\tilde{B}(t)=\Phi_{\tilde{N}(t)}\big((tR,\epsilon)\times\Sigma\big)$.
\begin{itemize}
\item[\emph{(a)}] There exists a diffeomorphism $\tilde{\upsilon}(t):\nu\big(\tilde{C}(t)\big)\rightarrow \nu\big(\tilde{B}(t)\big)$
which identifies the zero section in each bundle and maps $\nu_{r\sigma}\big(\tilde{C}(t)\big)$ to $\nu_{\Phi_{\tilde{N}(t)}(r,\sigma)}
\big(\tilde{B}(t)\big)$ for all $(r,\sigma)\in(tR,\epsilon)\times\Sigma$.
\item[\emph{(b)}] There exist an open subset $\tilde{V}(t)$ of the normal bundle
$\nu\big(tK_A\sqcup\tilde{B}(t)\big)$ of $tK_A\sqcup\tilde{B}(t)$ in $\R^7$, containing the zero section,
an open set $\tilde{S}(t)$ in $B(0;\epsilon_M)$, containing $tK_A\sqcup\tilde{B}(t)$, and a diffeomorphism $\tilde{n}(t)
:\tilde{V}(t)\rightarrow \tilde{S}(t)$.
Moreover, $\tilde{V}(t)$, $\tilde{S}(t)$ and $\tilde{n}(t)$ over $\tilde{B}(t)$ agree with $V_C$, $S_C$ and $n_C$ over $\tilde{C}(t)$, where
we use $\tilde{\upsilon}(t)$ in this identification, and $\tilde{V}(t)$ may be chosen to be an open neighbourhood 
of the zero section in $C^1$.
\end{itemize}
\end{cor}

\begin{proof}  Notice that, in the notation of Definition \ref{s4dfn2},
\begin{align*}
\Phi_{\tilde{N}(t)}(r,\sigma)-\iota(r,\sigma)&=t\big(1-f_{\text{inc}}(2t^{-\nu}r-1)\big)\big(\Phi_A(t^{-1}r,\sigma)-\iota(t^{-1}r,\sigma)\big)\\
&\;\,\;\,+f_{\text{inc}}(2t^{-\nu}r-1)\big(\Phi_N(r,\sigma)-\iota(r,\sigma)\big).
\end{align*}
The existence of 
$\tilde{\upsilon}(t)$ is then immediate from \eq{ch8s1eq1a} and \eq{dt1eq3a}. 
Let $V_{\tilde{C}(t)}=\{v\in \nu\big(\tilde{C}(t)\big)\,:\,v\in V_C\}$ 
and let $\tilde{S}^{\prime}(t)=n_C(V_{\tilde{C}(t)})$, which is an open subset of $\R^7$ containing $\tilde{C}(t)$.  
By making $V_C$ and $S_C$ smaller if necessary, but keeping the dilation invariance of the sets, we can ensure that $\tilde{S}^{\prime}(t)$ lies
within $B(0;\epsilon_M)$.  Moreover, in the notation of Corollaries \ref{s5cor1} and \ref{s5cor2},
$t^{-1}V_{\tilde{C}(t)}$ agrees with $V_A$ over $Q_A(t)=\Phi_A\big((R,t^{-1}\epsilon)\times\Sigma\big)$ via $\upsilon_A$, and
 $V_{\tilde{C}(t)}$ coincides with $V_B$ over 
$Q_N(t)=\Phi_N\big((tR,\epsilon)\times\Sigma\big)$ via $\upsilon_B$.  
Furthermore, $\tilde{S}^{\prime}(t)$ contains $tQ_A(t)$ and $Q_N(t)$.  It is therefore clear
that $\tilde{S}^{\prime}(t)$ contains $\tilde{B}(t)$.  Defining 
$\tilde{V}^{\prime}(t)=\tilde{\upsilon}(t)(V_{\tilde{C}(t)})$, it is also evident that
$\tilde{n}^{\prime}(t)=n_C\circ\tilde{\upsilon}(t)^{-1}$ 
is a diffeomorphism between $\tilde{V}^{\prime}(t)$ and $\tilde{S}^{\prime}(t)$.
We may then extend $\tilde{V}^{\prime}(t)$, $\tilde{S}^{\prime}(t)$ and $\tilde{n}^{\prime}(t)$ over $tK_A$ to the sets and map we require.
\end{proof}

We have thus constructed a tubular neighbourhood of $\chi^{-1}\big(\tilde{N}(t)\setminus K_N\big)$ that is conical over $\tilde{B}(t)$, 
which corresponds to the cone-like portion of $\tilde{N}(t)$.  So that we may relate this result to the normal bundle of $\tilde{N}(t)$, we make an
elementary observation.

\begin{lemma}\label{s5prop0} Use the notation of $\S$\ref{s4subs1}(b) and let $P$ be a 4-dimensional submanifold of $B(0;\epsilon_M)$.
There is a diffeomorphism $\chi_{\nu}$ from $\nu(P)$ to $\nu\big(\chi(P)\big)$, 
induced by $\chi$, 
such that $\chi_{\nu}|_p:\nu_p(P)\rightarrow\nu_{\chi(p)}\big(\chi(P)\big)$ for all $p\in P$.
\end{lemma}

\subsection{Neighbourhood of {\boldmath $\tilde{N}(t)$}}

Before our main result on tubular neighbourhoods, which follows from the work in $\S$\ref{s5subs1}, we make a couple of definitions.

\begin{dfn}\label{s5dfn1}
Let $P$ be a compact oriented 4-fold in $M$. 
Suppose further that 
$\big\|\varphi|_P\big\|_{C^0}<\epsilon_{\varphi}$ so that
Proposition \ref{s2prop4} holds for $P$, and
let $\jmath_P$
be the isomorphism between $\nu(P)$ and
$\Lambda^2_{\jmath+}T^*P=(\Lambda^2_+)_{g_P}T^*P$, using the
notation of that proposition. Suppose that $U_P$ is an open
set in $(\Lambda^2_+)_{g_P}T^*P$ containing the zero section, $T_P$
is an open set in $M$ containing $P$ and $\delta_P:U_P\rightarrow
T_P$ is a diffeomorphism such that $\delta_P(x,0)=x$ for all $x\in
P$.   We have a splitting $TU_P|_{(x,0)}=T_xP\oplus
(\Lambda^2_+)_{g_P}T^*_xP$ for all $x\in P$.  Thus we can consider
$d\delta_P$ at $P$ as a map from $TP\oplus(\Lambda^2_+)_{g_P}T^*P$
to $TP\oplus\nu(P)\cong TM|_{P}$. Hence, suppose further that, in
matrix notation,
\begin{equation*}\label{ch8s2subs1eq2}
d\delta_P|_{P}=\left(\begin{array}{cc}\mathfrak{I}& \mathfrak{A}\\
0 & \jmath_P^{-1}
\end{array}\right),\end{equation*}
where $\mathfrak{I}$ is the identity and $\mathfrak{A}$ is
arbitrary.  We say that such a diffeomorphism $\delta_P$ is
\emph{compatible with $\jmath_P$}.
\end{dfn}

\begin{note}
This compatibility is required so that we can use a result from \cite{McLean}, in $\S$\ref{s6}, 
which will help us describe the linearisation of the
deformation map at zero.
\end{note}

\begin{dfn}\label{TN} We remind the reader of the notation in $\S$\ref{s4subs1}(b).
By \cite[Proposition 6.4]{Lotay2}, there exist a tubular neighbourhood $\hat{T}$ of $\hat{N}$
in $M$, an open set $\hat{U}\subseteq(\Lambda^2_+)_{g|_{\hat{N}}}T^*\hat{N}$ and a diffeomorphism $\hat{\delta}:
\hat{U}\rightarrow\hat{T}$ compatible with $\jmath_{\hat{N}}$, in the notation of Proposition \ref{s2prop2} and Definition \ref{s5dfn1}.
Define $$T_N=\hat{T}\cup V$$ and note that it is an open subset of $M$ containing $N$.  
Moreover, by \cite[Proposition 6.19]{Lotay2}, $T_N$ retracts onto $N$ and $H^3_{\text{dR}}(T_N)$ is isomorphic to
$H^3_{\text{cs}}(\hat{N})$.  As a consequence of this, 
$\varphi|_{T_N}$ is exact. \end{dfn}

Recall, by Lemma \ref{smallphi},
 that if $\tau$ is small enough, $\tilde{N}(t)$ will satisfy the conditions of Proposition \ref{s2prop4} for all $t\in(0,\tau)$.

\begin{prop}\label{s5prop3} Use the notation of $\S$\ref{s4subs1}-\ref{s4subs2}, 
Corollary \ref{s5cor4} and Lemma \ref{s5prop0}.
Let $\tilde{\jmath}(t)=\jmath_{\tilde{N}(t)}$ and $\tilde{g}(t)=g_{\tilde{N}(t)}$, in the notation of Proposition \ref{s2prop4}, and
let $\rho_t$ be a radius function on $\tilde{N}(t)$ as given in Definition \ref{s4dfn3}.
There exist an open subset $\tilde{U}(t)$
of $(\Lambda^2_+)_{\tilde{g}(t)}T^*\tilde{N}(t)$ containing the zero section and 
$\tilde{W}(t)=\big(\tilde{\jmath}(t)\circ\chi_{\nu}\big)\big(\tilde{V}(t)\big)$, an open subset $\tilde{T}(t)$ of
$M$ containing $\tilde{N}(t)$ and $\chi\big(\tilde{S}(t)\big)$, and a diffeomorphism $\tilde{\delta}(t):\tilde{U}(t)\rightarrow\tilde{T}(t)$
which is compatible with $\tilde{\jmath}(t)$, in the sense of Definition \ref{s5dfn1}, and such that the following diagram commutes:
\begin{equation}\label{s5eq1}
\begin{gathered}
\xymatrix{ \tilde{W}(t)\ar[rr]^{\tilde{\jmath}(t)^{-1}}\ar[dd]_{\tilde{\delta}(t)}
&&\chi_{\nu}\big(\tilde{V}(t)\big) \ar[d]^{\chi_{\nu}^{-1}}
\\
&&\tilde{V}(t)\ar[d]^{\tilde{n}(t)} \\
 \chi\big(\tilde{S}(t)\big) & &\tilde{S}(t).\ar[ll]_{\chi}
}\end{gathered}
\end{equation}
Moreover, $\tilde{U}(t)$ can be chosen to be an open neighbourhood of the zero section in $C^1_{1,\,t}$, given in Definition \ref{s4dfn4},
and $\tilde{T}(t)$ is contained in $T_N$, defined in Definition \ref{TN}.
\end{prop}

\begin{proof} We first see, by Corollary \ref{s5cor4} and Lemma \ref{s5prop0}, that 
\begin{align*}
\chi_{\nu}\big(\tilde{V}(t)\big)\subseteq\nu\Big(\chi\big(tK_A\sqcup\tilde{B}(t)\big)\Big)
=\nu\big(\tilde{N}(t)\setminus K_N\big)
\end{align*}
and is a $C^1$-open neighbourhood of the zero section over $\tilde{N}(t)\setminus K_N$.  So $\tilde{W}(t)\subseteq(\Lambda^2_+)_{\tilde{g}(t)}
T^*\tilde{N}(t)$ and we may
define $\tilde{\delta}(t)|_{\tilde{W}(t)}$ by the commutative diagram \eq{s5eq1}.  Note that $\tilde{\delta}(t)$ acts as the identity
on the zero section and, by the choice of $\chi_{\nu}$ in Lemma \ref{s5prop0}, 
satisfies the compatibility condition with $\tilde{\jmath}(t)$ over $\tilde{N}(t)\setminus K_N$. We may thus extend $\tilde{W}(t)$,
$\chi\big(\tilde{S}(t)\big)$ and $\tilde{\delta}(t)|_{\tilde{W}(t)}$ smoothly over the compact set $K_N$
to the open sets $\tilde{U}(t)$ and $\tilde{T}(t)$ and the diffeomorphism $\tilde{\delta}(t)$ we need.  
Notice that $\tilde{S}(t)\subseteq B(0;\epsilon_M)$ implies that $\chi\big(\tilde{S}(t)\big)\subseteq V\subseteq T_N$.  We may easily ensure 
that $\tilde{U}(t)$ is an open neighbourhood of the zero section in $C^1_{1,\,t}$ by the construction of $\tilde{V}(t)$
in the proof of Corollary \ref{s5cor4}.
\end{proof}

\begin{note}
We shall use the notation of Definition \ref{TN} and Proposition \ref{s5prop3} throughout the sequel, though we shall make 
every attempt to notify the reader when we do this.
\end{note}

\subsection{Bounds and estimates}

In this subsection, we take a first look at the behaviour of various geometric objects on $\tilde{N}(t)$ as $t$ varies.  We start with
the size of elements of $\tilde{U}(t)$.

\begin{dfn}\label{s5dfn3} Use the notation of Definitions \ref{s4dfn2}-\ref{s4dfn4} and Proposition \ref{s5prop3}.
Since $\tilde{U}(t)$ is an open neighbourhood of the zero section in $C^1_{1,\,t}$, 
it is an open set in $L^p_{2,\,1,\,t}$ for $p>4$. 
Therefore, by note (b) after Definition \ref{s4dfn4}, there exists a constant $\tilde{\epsilon}>0$ such that
if $\alpha\in L^8_{2,\,1,\,t}\big((\Lambda^2_+)_{\tilde{g}(t)}T^*\tilde{N}(t)\big)$
 with $\|\alpha\|_{L^8_{2,\,1,\,t}}\leq \tilde{\epsilon}$, 
$\alpha\in L^8_{2,\,1,\,t}\big(\tilde{U}(t)\big)$.  Moreover, by the $t$ dependence of the norm on $L^8_{2,\,1,\,t}$ and
the construction of $\tilde{U}(t)$, we can ensure that $\tilde{\epsilon}$ is independent of $t$.  
In particular, notice that $\|\alpha\|_{L^8_{2,\,1,\,t}}$ being small ensures that $\|\rho_t^{-1}\alpha\|_{C^0}$
and $\|\nabla_t\alpha\|_{C^0}$ are small.
\end{dfn}

We make another definition.

\begin{dfn}\label{s5dfn4} Use the notation of Proposition \ref{s5prop3}.
Define $\tilde{\varphi}(t)=\tilde{\delta}(t)^*(\varphi)$, which is a 3-form on $\tilde{U}(t)$.
\end{dfn}

\noindent We can estimate the size of $\tilde{\varphi}(t)$ and its derivatives as follows.

\begin{prop}\label{s5prop4} Suppose that $\tau$ is sufficiently small.  Use the notation of $\S$\ref{s4subs1}(b), Proposition \ref{s5prop3}
and Definition \ref{s5dfn4}.
\begin{itemize}
\item[\emph{(a)}] For all $j\in\N$ there exists a constant $C(\tilde{\varphi})_j>0$, independent of $t$,
such that $$|\nabla^j\tilde{\varphi}(t)|\leq
C(\tilde{\varphi})_j\rho_t^{-j},$$ where $\nabla$ is the Levi--Civita connection of $g$, the metric on $M$.
\item[\emph{(b)}] The injectivity radius $\delta\big(\tilde{g}(t)\big)$
and Riemann curvature $\text{\emph{Riem}}\big(\tilde{g}(t)\big)$ of 
$\tilde{g}(t)$ on $\tilde{N}(t)$ satisfy $$\big|\,\delta\big(\tilde{g}(t)\big)\big|\geq
C\big(\delta(\tilde{g})\big)\rho_t\quad\text{and}\quad\big|\,\text{\emph{Riem}}\big(\tilde{g}(t)\big)\big|\leq
C\big(\text{\emph{Riem}}(\tilde{g})\big)\rho_t^{-2}$$ for some constants
$C\big(\delta(\tilde{g})\big)>0$ and $C\big(\text{\emph{Riem}}(\tilde{g})\big)>0$ independent of $t$.
\end{itemize}
\end{prop}

\begin{proof}  Recall the notation of $\S$\ref{s4subs1}(a)-(c) and Definition \ref{s4dfn2} and let 
$g_0$ be the Euclidean metric on $\R^7$.
 The dominant contributions to each of the quantities
in which we are interested comes from its behaviour near
$\chi(tK_A)$, as long as $t$ is
sufficiently small. However, we know that here 
$\tilde{g}(t)\approx t^2g_0$, whereas the magnitude of the metric over $\Psi_{\tilde{N}(t)}\big((tR,\epsilon)\times\Sigma\big)$
depends on $r$.  Moreover, we see that near $\chi(tK_A)$, $\tilde{\varphi}(t)\approx
t^3\varphi_0$, $\delta\big(\tilde{g}(t)\big)\approx
t\delta(g_0)$ and $\text{Riem}\big(\tilde{g}(t)\big)\approx
t^2\text{Riem}(g_0)$ as long as $\tau$ is sufficiently
small.  Similar estimates hold involving $r$ instead of $t$ over $\Psi_{\tilde{N}(t)}\big((tR,\epsilon)\times\Sigma\big)$.
The result follows by the definition of $\rho_t$.
\end{proof}

We now turn to interior regularity estimates for the operator $d+d^*$ acting on self-dual 2-forms and 4-forms.
The first result is a direct corollary of \cite[Theorem 6.2.6]{Morrey}, which is a result for linear elliptic systems 
for balls in Euclidean space.

\begin{prop}\label{s5prop6}  Let $p>1$ and let $k\in\N$.
Let $B_s$ denote the ball radius $s>0$ in $\R^4$ and let $\nabla$ denote the
Levi--Connection of the Euclidean metric on $\R^4$.
There exist constants $\epsilon_{\text{\emph{reg}}}>0$ and $C(d+d^*)_0>0$ 
such that if $\kappa\in(0,1]$ and $(\alpha,\beta)\in L^p_{k+1}\big(\Lambda^2_+T^*B_{3\kappa\epsilon_{\text{\emph{reg}}}}
\oplus\Lambda^4T^*B_{3\kappa\epsilon_{\text{\emph{reg}}}}\big)$ ,
\begin{align}
\sum_{j=0}^{k+1}&\kappa^{j-\frac{4}{p}}\big\|\nabla^j(\alpha,\beta)|_{B_{2\kappa\epsilon_{\text{\emph{reg}}}}}\big\|_{L^p}
\nonumber\\
&\leq C(d+d^*)_0\left(\sum_{j=0}^k\kappa^{j+1-\frac{4}{p}}\|\nabla^j(d\alpha+d^*\beta)\|_{L^p}
+\kappa^{-4}\|(\alpha,\beta)\|_{L^1}\right).\label{s5prop8eq0}
\end{align}
\end{prop}

\begin{proof}
The linear first-order differential operator $\mathcal{L}:\gamma=(\alpha,\beta)\mapsto d\alpha+d^*\beta$ is smooth and elliptic so, 
by \cite[Theorem 6.2.6]{Morrey}, there exist constants $\epsilon_{\text{reg}}>0$ and $C(\mathcal{L})>0$ such that, for
all $0<s\leq \epsilon_{\text{reg}}$ and $\gamma\in L^p_{k+1}(\Lambda^2_+T^*B_{3s}\oplus\Lambda^4T^*B_{3s})$, 
\begin{align*}
\sum_{j=0}^{k+1}(3s)^{j-k-1}&\|\nabla^j\gamma|_{B_{2s}}\|_{L^p}\\
&\leq C(\mathcal{L})\Bigg(
\sum_{j=0}^k(3s)^{j-k}\|\nabla^j(\mathcal{L}\gamma)\|_{L^p}
+(3s)^{-k-1-4\left(\frac{p-1}{p}\right)}\|\gamma\|_{L^1}\Bigg).
\end{align*}
Multiplying through by $(3s)^{k+1-\frac{4}{p}}$ and substituting $s=\kappa\epsilon_{\text{reg}}$, we have the result 
for a constant $C(d+d^*)_0$ depending on $C(\mathcal{L})$ and $\epsilon_{\text{reg}}$.
\end{proof}


\noindent Our next lemma is necessary since we need to bound variations in the radius function $\rho_t$ in 
small balls on $\tilde{N}(t)$.

\begin{lemma}\label{lip}  Use the notation of Definitions \ref{s4dfn2}, \ref{s4dfn3} and \ref{s4dfn4}.
For $x\in\tilde{N}(t)$ and $s>0$, denote by
$B_s(x)$ the geodesic ball of radius $s$ in $\tilde{N}(t)$ with respect to $\tilde{g}(t)$.  There exists a constant 
$c_0>0$, independent of $t$, such that for all $x\in\tilde{N}(t)$ and $y\in B_{c_0\rho_t(x)}(x)$,
$$|\rho_t(y)-\rho_t(x)|\leq \frac{1}{2}\rho_t(x).$$
\end{lemma}

\begin{proof}
The function $\rho_t$ is smooth on the compact manifold $\tilde{N}(t)$ and so is Lipschitz.  Therefore, 
if $d(x,y)$ is the geodesic distance between $x,y\in\tilde{N}(t)$, there exists $a>0$ such
that
$$|\rho_t(x)-\rho_t(y)|\leq a\, d(x,y)$$
for all $x,y\in\tilde{N}(t)$.  As can be seen from the definition of $\rho_t$ in Definition \ref{s4dfn3}, 
there is a $t$ independent bound for the modulus of the derivative of $\rho_t$.
Thus, $a$ can be chosen to be independent of $t$ and we take $c_0>0$ to be such that $ac_0\leq\frac{1}{2}$.
\end{proof}

\noindent The next result shows that in small balls on $\tilde{N}(t)$ the metric $\tilde{g}(t)$ is ``close" to the
Euclidean metric.  

\begin{prop}\label{s5prop7}
Use the notation of Definitions \ref{s4dfn2}, \ref{s4dfn3} and \ref{s4dfn4}, 
Propositions \ref{s5prop4} and \ref{s5prop6} and Lemma \ref{lip}.   
Let $c_1>0$ be a sufficiently small constant.  
There exist constants $c_2\in(0,1)$, $c_3>0$ and $c_4>0$, depending only on $c_0$, $c_1$, $C\big(\delta(\tilde{g})\big)$ and
 $C\big(\text{\emph{Riem}}(\tilde{g})\big)$, such that the following holds.  
For $x\in\tilde{N}(t)$ let $\eta_t(x)=c_2\epsilon_{\text{\emph{reg}}}\rho_t(x)$ and let 
$4c_2\epsilon_{\text{\emph{reg}}}\leq c_0$.  There exist smooth injective maps
$\Upsilon_x:B_{3\eta_t(x)}\rightarrow \tilde{N}(t)$ such that 
\begin{gather*}
\big\|
\Upsilon_x^*\big(\tilde{g}(t)\big)-g_0\big\|_{L^8_2}\leq c_1,\\
B_{\eta_t(x)}(x)\subseteq \Upsilon_x\big(B_{2\eta_t(x)}\big)\subseteq \Upsilon_x\big(B_{3\eta_t(x)}\big)\subseteq B_{4\eta_t(x)}(x)\;\,\text{and}\\
c_3\{\rho_t(x)\}^4\leq \text{\emph{vol}}\Big(\Upsilon_x\big(B_{2\eta_t(x)}\big)\Big)\leq 
\text{\emph{vol}}\Big(\Upsilon_x\big(B_{3\eta_t(x)}\big)\Big)
\leq c_4\{\rho_t(x)\}^4.
\end{gather*}
\end{prop}

\begin{proof}
This is an analogous result to \cite[Proposition 5.9]{Joyce3} and can be proved in an identical manner.  Therefore we only 
sketch the proof here.  By Proposition \ref{s5prop4}(b),
 the injective radius is bounded below and the sectional curvature is bounded above on $(\tilde{N}(t),\tilde{g}(t))$.
We can then use Jost and Karcher's \cite{Jost} theory of harmonic coordinates on $\tilde{N}(t)$, which give a
$C^{1,\,a}$ estimate for the metric,
and the improvement of these results to $L^q_2$ for $q>2$ by Anderson, which are described in Petersen \cite[$\S$4 \& $\S$5]{Petersen}.  
These give our coordinate systems $\Upsilon_x$
and our first estimate.  Notice that, unlike in \cite[Proposition 5.9]{Joyce3}, 
we do not need to rescale the metric $\Upsilon_x^*\big(\tilde{g}(t)\big)$ by 
$\{\eta_t(x)\}^{-2}$ to make it close to $g_0$ on $B_{3\eta_t(x)}$, since here the balls in Euclidean space scale with $\eta_t(x)$.
The volume estimates and the nesting of the balls in $\tilde{N}(t)$ follow from the fact that an $L^8_2$
 estimate on the metric ensures that we have control on the $C^0$ norm of the metric.
%
\end{proof}

\noindent We can now prove our interior regularity estimate in a similar manner to \cite[Proposition 5.11]{Joyce3}.
The proposition is an analogue of the elliptic regularity result for weighted Banach spaces given
in Proposition \ref{s3sub4prop3}.

\begin{prop}\label{s5prop8}  Use the notation of Definitions \ref{s4dfn2} and \ref{s4dfn4}.  Let $1<p\leq 8$, let
$k=0$ or $1$ and let $\eta\in\R$.
There exists a constant $C(d+d^*)>0$, independent of $t$, such that, for all $\alpha\in L^p_{k+1}\big((\Lambda^2_+)_{\tilde{g}(t)}
T^*\tilde{N}(t)\big)$,
$$\|\alpha\|_{L^p_{k+1,\,\eta,\,t}}\leq C(d+d^*)\Big(\|d\alpha\|_{L^p_{k,\,\eta-1,\,t}}+\|\alpha\|_{L^1_{0,\,\eta,\,t}}\Big)$$
\end{prop}

\begin{proof} 
Let $c_1>0$ be a sufficiently small constant, apply Proposition \ref{s5prop7} and let $x\in\tilde{N}(t)$.  
Use the notation of that result along with
the notation of Propositions \ref{s5prop6} and Lemma \ref{lip}.  Let $g_x=\Upsilon_x^*\big(\tilde{g}(t)\big)$,
which is a metric on $B_{3\eta_t(x)}$.  

For 
$\alpha\in L^p_{k+1}\big((\Lambda^2_+)_{\tilde{g}(t)}
T^*\tilde{N}(t)\big)$, define 
$$\gamma_0=(\Upsilon_x^*(\alpha),0)\in L^p_{k+1}\big((\Lambda^2_+)_{g_x}T^*B_{3\eta_t(x)}
\oplus \Lambda^4T^*B_{3\eta_t(x)}\big).$$
Recall that $\eta_t(x)=c_2\rho_t(x)\epsilon_{\text{reg}}$ 
and $\kappa=c_2\rho_t(x)\in(0,1)$ since $c_2\in(0,1)$ and $\rho_t(x)\in (0,1]$.

If $c_1$ is sufficiently small, $g_x$ and $g_0$ are close in $L^8_2$ by Proposition \ref{s5prop7}, 
thus close in $L^p_{k+1}$ by the Sobolev Embedding Theorem.  
Hence we can increase $C(d+d^*)_0$ to
$C(d+d^*)_1$, depending only on $c_1$, $C(d+d^*)_0$ and $\epsilon_{\text{reg}}$, 
such that Proposition \ref{s5prop6} holds for $(\alpha_0,\beta_0)\in
L^p_{k+1}\big((\Lambda^2_+)_{g_x}T^*B_{3\eta_t(x)}\oplus \Lambda^4T^*B_{3\eta_t(x)}\big)$, with moduli, Hodge star and
covariant derivatives calculated using the metric $g_x$.  Note that we are taking $\kappa=c_2\rho_t(x)$ in 
\eq{s5prop8eq0}.  
Putting $\gamma_0$ into \eq{s5prop8eq0}, pushing
 forward by $\Upsilon_x$ and multiplying by $\{\rho_t(x)\}^{-\eta}$ gives that
\begin{align}
\sum_{j=0}^{k+1}
&\big\|\{\rho_t(x)\}^{j-\eta-\frac{4}{p}}\nabla_t^j\alpha|_{\Upsilon_x\big(B_{2\eta_t(x)}\big)}\big\|_{L^p}
\nonumber\\
&\leq C(d+d^*)_2 \Bigg(\sum_{j=0}^k
\big\|\{\rho_t(x)\}^{j+1-\eta-\frac{4}{p}}\nabla_t^j(d\alpha)|_{\Upsilon_x\big(B_{3\eta_t(x)}\big)}\big\|_{L^p}
\nonumber \\
&\qquad\qquad\qquad\qquad\qquad\qquad\qquad+\big\|\{\rho_t(x)\}^{-\eta-4}\alpha|_{\Upsilon_x\big(B_{3\eta_t(x)}\big)}
\big\|_{L^1}\Bigg), 
\label{s5prop8eq2}
\end{align}
where $C(d+d^*)_2>0$ is a constant which depends only on $C(d+d^*)_1$ and $c_2$, so independent of $x$ and $t$. 

Now, by Lemma \ref{lip}, $|\frac{\rho_t(y)}{\rho_t(x)}-1|\leq\frac{1}{2}$ for all 
$y\in B_{4\eta_t(x)}(x)$ since $4c_2\epsilon_{\text{reg}}\leq c_0$ by the choice of $c_2$ in
Proposition \ref{s5prop7}.  Moreover, by Proposition \ref{s5prop7}, $\Upsilon_x\big(B_{3\eta_t(x)}\big)
\subseteq B_{4\eta_t(x)}(x)$.  Notice that 
if we 
replace 
$\rho_t(x)$ by $\rho_t|_{\Upsilon_x\big(B_{m\eta_t(x)}\big)}$ ($m=2$ or $3$) in \eq{s5prop8eq2}, 
we get the weighted $L^p$ space norms as given in Definition \ref{s4dfn4} that we require, since
we have factors of $\rho_t^{-\frac{4}{p}}$ and $\rho_t^{-4}$ for the $L^p$ and $L^1$ norms respectively.  
In conclusion, we may increase the constant $C(d+d^*)_2$ to another constant $C(d+d^*)_3$, independent of $x$ and $t$,
 such that
\begin{align}
\big\|\alpha &|_{\Upsilon_x\big(B_{2\eta_t(x)}\big)}\big\|_{L^p_{k+1,\,\eta,\,t}}\nonumber\\
&\qquad\leq 
C(d+d^*)_3 \Big(\big\|d\alpha|_{\Upsilon_x\big(B_{3\eta_t(x)}\big)}\big\|_{L^p_{k,\,\eta-1,\,t}}
+\big\|\alpha|_{\Upsilon_x\big(B_{2\eta_t(x)}\big)}\big\|_{L^1_{0,\,\eta,\,t}}\Big).\label{s5prop8eq3}
\end{align}

To complete the proof we first use the compactness of $\tilde{N}(t)$ to finitely cover it using open sets of the form
$B_{\eta_t(x)}(x)$.  Recall the nesting of sets given in Proposition \ref{s5prop7} and the fact that $C(d+d^*)_3$
is independent of $x$ and $t$.  We can therefore use a partition of unity for the 
finite open cover by small geodesic balls together with \eq{s5prop8eq3} to provide our required estimate.
\end{proof}




\section{The deformation map}\label{s6}

The next stage is to consider \emph{deformations} of $\tilde{N}(t)$ which are \emph{coassociative}.
The key is to define a \emph{deformation map}.   We introduce the following notation: 
$$L^p_k\big(\tilde{U}(t)\big)=\big\{\xi\in L^p_k\big(\Lambda^mT^*\tilde{N}(t)\big)\,:\,G_{\xi}\subseteq \tilde{U}(t)\big\},$$
where $G_\xi$ denotes the graph of the form $\xi$.  Note that for this definition to make sense we require that the forms
be continuous, so that their graphs are well-defined at all points; that is, we need $k>\frac{4}{p}$.  We adopt similar notation
for subsets of $C^k$ and H\"older spaces.

\begin{dfn}\label{s6dfn1}  Use the notation of $\S$\ref{s4subs1}, Definition \ref{s4dfn2} and Proposition \ref{s5prop3}.
 Let $\alpha\in C^1\big(\tilde{U}(t)\big)$ and let $\pi_\alpha(t):\tilde{N}(t)\rightarrow G_\alpha$ 
 be given by $\pi_\alpha(t)(x)=(x,\alpha(x))$.  Define $f_\alpha(t)=\tilde{\delta}(t)\circ\pi_\alpha(t)$ and let
$\tilde{N}_\alpha(t) =f_\alpha(t)\big(\tilde{N}(t)\big)\subseteq\tilde{T}(t)$.  The \emph{deformation map},
$F_t:C^1\big(\tilde{U}(t)\big)\rightarrow
C^0\big(\Lambda^3T^*\tilde{N}(t)\big)$, is given by
$$F_t(\alpha)=f_\alpha(t)^*\left(\varphi|_{\tilde{N}_\alpha(t)}\right).$$
Immediately note that $F_t(0)=\varphi|_{\tilde{N}(t)}$ and that $F_t$ is a
smooth map. Furthermore, by \cite[p. 731]{McLean}, which we are
allowed to use by virtue of the compatibility of $\tilde{\delta}(t)$ with $\tilde{\jmath}(t)$ and the
fact that $\varphi$ is closed, the linearisation of $F_t$ at $0$ is
$$dF_t|_0(\alpha)=d\alpha$$
for $\alpha\in C^1\big((\Lambda^2_{+})_{\tilde{g}(t)}T^*\tilde{N}(t)\big)$.

It is clear that $\text{Ker}\,F_t$ is equal to the set of $\alpha\in C^1\big(\tilde{U}(t)\big)$ such that $\tilde{N}_\alpha(t)$
is coassociative.
\end{dfn}

The first thing we must do is rewrite $F_t$.

\begin{prop}\label{s6prop1}  Use the notation of Definitions \ref{s4dfn2}-\ref{s4dfn4}, 
\ref{s5dfn3} and \ref{s6dfn1}
 and Proposition \ref{s5prop3}.
The map $F_t$ can be
written as
\begin{equation}\label{PF}
F_t(\alpha)(x)=\varphi(x)+d\alpha(x) +P_{F_t}\big(x,\alpha(x),\nabla_t\alpha(x)\big)
\end{equation}
for $x\in \tilde{N}(t)$, where $$P_{F_t}:\{(x,y,z)\,:\,(x,y)\in
\tilde{U}(t),\,z\in
T_x^*\tilde{N}(t)\otimes(\Lambda_{+}^2)_{\tilde{g}(t)}T_x^*\tilde{N}(t)\}\rightarrow\Lambda^3T^*\tilde{N}(t)$$
is a smooth map such that $P_{F_t}(x,y,z)\in\Lambda^3T^*_x\tilde{N}(t)$.
Denote
$P_{F_t}\big(x,\alpha(x),\nabla\alpha(x)\big)$ by $P_{F_t}(\alpha)(x)$ for all $x\in\tilde{N}(t)$.
For $p\geq 1$ and $k\in\N$, if $\alpha\in L^{p}_{k+1}\big(\tilde{U}(t)\big)$ with
$\|\alpha\|_{C^1_{1,\,t}}$ sufficiently small, $P_{F_t}(\alpha)\in L^{p}_k\big(\Lambda^3T^*\tilde{N}(t)\big)$.  
Similarly, for $l\in\N$ and $a\in(0,1)$, if $\alpha\in C^{l+1,\,a}\big(\tilde{U}(t)\big)$ with $\|\alpha\|_{C^1_{1,\,t}}$ sufficiently
small, $P_{F_t}(\alpha)\in C^{l,\,a}\big(\Lambda^3T^*\tilde{N}(t)\big)$.  

Moreover, suppose $\alpha\in L^8_{2,\,1,\,t}\big(\tilde{U}(t)\big)$ with $\|\alpha\|_{L^8_{2,\,1,\,t}}\leq \tilde{\epsilon}$.
Making $\tilde{\epsilon}$ smaller if necessary, there exist constants $p_0>0$ and $p_1>0$, independent of $\alpha$ and $t$
such that
\begin{align*}
\|P_{F_t}(\alpha)\|_{L^{4/3}}\leq p_0\|\alpha\|_{L^{4/3}_{1,\,-2,\,t}}^2
\quad\text{and}\quad\|P_{F_t}(\alpha)\|_{L^8_{1,\,0,\,t}}\leq p_1\|\alpha\|_{L^8_{2,\,1,\,t}}^2.
\end{align*}
\end{prop}

\begin{proof} Recall the notation of $\S$\ref{s4subs1}-\ref{s4subs2}.  Note that, if $G_{\alpha}$ is the graph of $\alpha$, 
$F_t(\alpha)(x)$ depends on
$T_{(x,\alpha(x))}G_{\alpha}$,  hence on $\alpha(x)$ and $\nabla_t\alpha(x)$.  We may
thus define the map $P_{F_t}$ by \eq{PF}.
The rest of this proposition is proved in a similar manner to \cite[Proposition
4.3]{Lotay1} and \cite[Proposition 6.9]{Lotay2}, except here the
analysis is more straightforward since $\tilde{N}(t)$ is compact.  The proof rests on the
observation that $P_{F_t}$ and its first derivatives are zero at
$\alpha=0$ and thus it is approximately quadratic in $\alpha$ and
$\nabla_t\alpha$, as long as $\rho_t^{-1}\alpha$ and $\nabla_t\alpha$ are small in $C^0$; i.e. $\|\alpha\|_{C^1_{1,\,t}}$ 
is small.  The power of $\rho_t$
is determined so that $P_{F_t}(\alpha)$ scales as expected under variations in $t$.  Therefore, if $\alpha\in L^p_{k+1}$ or $C^{l+1,\,a}$
and $\|\alpha\|_{C^1_{1,\,t}}$ is sufficiently small, $P_{F_t}(\alpha)$ lies in $L^p_k$ or $C^{l,\,a}$ respectively.

Furthermore, if $\alpha$ has at least one derivative and is near zero in $C^1_{1,\,t}$, 
\begin{align}
|P_{F_t}(\alpha)|&=O\Big(\big(|\rho_t^{-1}\alpha|+|\nabla_t\alpha|\big)^2\Big)\label{s6prop1eq1}\\
\intertext{and, similarly, if $\alpha$ has two derivatives,}
\big|\nabla_t\big(P_{F_t}(\alpha)\big)\big|&=O\Big(\rho_t^{-1}\big(|\rho_t^{-1}\alpha|+|\nabla_t\alpha|+|\rho_t\nabla_t^2\alpha|\big)^2\Big).\label{s6prop1eq2}
\end{align}
There is an obvious generalisation to higher derivatives, but we shall not require this.  Integrating \eq{s6prop1eq1}
we see that
$$\|P_{F_t}(\alpha)\|_{L^{4/3}}=O\Big(\big(\|\rho_t^{-1}\alpha\|_{L^{4/3}}+\|\nabla_t\alpha\|_{L^{4/3}}
\big)^2\Big),$$
which gives our first estimate, recalling the definition of the weighted norm in Definition \ref{s4dfn4}.  Using \eq{s6prop1eq1}
and \eq{s6prop1eq2}, we have the stronger result
$$\|P_{F_t}(\alpha)\|_{L^8_{1,\,0,\,t}}=O\Big(\|\alpha\|^2_{L^8_{2,\,\frac{3}{4},\,t}}\Big),$$
where $\frac{3}{4}$ is calculated as $1-\frac{2}{8}$.  Using the fact that $|\rho_t|\leq 1$, our second estimate follows.
\end{proof}

We now make a further estimate which builds upon the previous proposition.

\begin{prop}\label{s6prop3}  Use the notation of Definitions \ref{s4dfn4} and \ref{s5dfn3} and Propositions \ref{s5prop3} and \ref{s6prop1}.
Let $\alpha,\beta\in L^8_{2,\,1,\,t}\big(\tilde{U}(t)\big)$
and suppose that 
$\|\alpha\|_{L^8_{2,\,1,\,t}},\|\beta\|_{L^8_{2,\,1,\,t}}\leq \tilde{\epsilon}$.  Making $\tilde{\epsilon}$ smaller if
necessary, there exists a constant $C(P_{F_t})>0$, independent of $\alpha$, $\beta$ and $t$, such that
\begin{align*}
\|P_{F_t}(\alpha)-P_{F_t}(\beta)\|_{L^{4/3}}&\leq C(P_{F_t})\|\alpha-\beta\|_{L^{4/3}_{1,\,-2,\,t}}
\bigg(\|\alpha\|_{L^{4/3}_{1,\,-2,\,t}}+\|\beta\|_{L^{4/3}_{1,\,-2,\,t}}\bigg) \\
\intertext{and}
\|P_{F_t}(\alpha)-P_{F_t}(\beta)\|_{L^8_{1,\,0,\,t}}&\leq C(P_{F_t})\|\alpha-\beta\|_{L^8_{2,\,1,\,t}}
\Big(\|\alpha\|_{L^8_{2,\,1,\,t}}+\|\beta\|_{L^8_{2,\,1,\,t}}\Big)
.
\end{align*}
\end{prop}

\noindent This result is analogous to \cite[Proposition 5.8]{Joyce3} and the proof is almost identical.  Rather than repeating that technical proof with only cosmetic changes, we hope to convince the reader with a heuristic argument.  If $q$ is a homogeneous quadratic, 
it is evident that there exists a constant $C(q)>0$ such that
$$|q(x)-q(y)|\leq C(q)|x-y|(|x|+|y|)$$
for all $x$ and $y$.  In the proof of our previous proposition we showed that $P_{F_t}$ is approximately quadratic for forms which are sufficiently
 near zero in $C^1_{1,\,t}$, which is ensured by the forms being small in the $L^8_{2,\,1,\,t}$ norm.  The introduction of powers of $\rho_t^{-1}$ in the estimates above is necessary since the bounds on moduli of certain
 derivatives of $F_t$ are dependent upon the moduli of corresponding derivatives, up to a maximum of third order, of $\tilde{\varphi}(t)$,
as given in Definition \ref{s5dfn4}.  These
 factors are thus determined by Proposition \ref{s5prop4}(a).

We can now provide further useful estimates by considering a map associated to $F_t$.  

\begin{dfn}\label{s6dfn5} Use the notation of Definitions \ref{s4dfn2}-\ref{s4dfn4}, \ref{TN} and \ref{s6dfn1},
Corollary \ref{s2cor1} 
and Proposition \ref{s5prop3}.  Let $p\geq 1$ and let $k\in\N$.  Recall that, for $\alpha\in L^p_{k+1}\big(\tilde{U}(t)\big)$,
 $N_{\alpha}(t)\subseteq T_N$ and $\varphi|_{T_N}$ is exact.  Therefore, if $\|\alpha\|_{C^1_{1,\,t}}$ is sufficiently small, 
we may apply Proposition \ref{s6prop1} and deduce that
$F_t(\alpha)\in L^p_k\big(\Lambda^3T^*\tilde{N}(t)\big)$ and is exact.  By Corollary \ref{s2cor1}, we may thus
define $$H_t(\alpha)=d^{-1}\big(F_t(\alpha)\big)\in V^p_{k+1}.$$
We see that $H_t$ is a smooth map with $d\big(H_t(\alpha)\big)=F_t(\alpha)$, $H_t(0)=\psi(t)\in V^p_{k+1}$ such
that $d\psi(t)=\varphi|_{\tilde{N}(t)}$ and $dH_t|_0(\alpha)=\alpha_{H_t}$, where $\alpha_{H_t}$ is
the $L^2$-orthogonal projection of $\alpha\in L^p_{k+1}\big((\Lambda^2_+)_{\tilde{g}(t)}T^*\tilde{N}(t)\big)$
onto $V^p_{k+1}$.
\end{dfn}

We have the analogue of Proposition \ref{s6prop1} for $H_t$.

\begin{prop}\label{s6prop2} Use the notation of Definitions \ref{s4dfn2}-\ref{s4dfn4} and \ref{s6dfn5},
Corollary \ref{s2cor1} and Proposition \ref{s5prop3}.  Let $p\geq 1$ and $k\in\N$.
The map $H_t$, for $\alpha\in L^p_{k+1}\big(\tilde{U}(t)\big)$ with $\|\alpha\|_{C^1_{1,\,t}}$ sufficiently small, can be
written as
\begin{equation}\label{PH}
H_t(\alpha)(x)=\psi(t)(x)+\alpha_{H_t}(x) +P_{H_t}\big(x,\alpha(x),\nabla_t\alpha(x)\big)
\end{equation}
for $x\in \tilde{N}(t)$, where $$P_{H_t}\!:\{(x,y,z):(x,y)\in
\tilde{U}(t),\,z\in
T_x^*\tilde{N}(t)\otimes(\Lambda_{+}^2)_{\tilde{g}(t)}T_x^*\tilde{N}(t)\}\rightarrow(\Lambda^2_+)_{\tilde{g}(t)}T^*\tilde{N}(t)$$
is a smooth map such that $P_{H_t}(x,y,z)\in(\Lambda^2_+)_{\tilde{g}(t)}T^*_x\tilde{N}(t)$.
Let $P_{H_t}(\alpha)(x)$ denote
$P_{H_t}\big(x,\alpha(x),\nabla_t\alpha(x)\big)$  for all $x\in\tilde{N}(t)$.    
Then
$P_{H_t}(\alpha)\in
V^p_{k+1}$.
\end{prop}

\begin{proof}
We can define the smooth map $P_{H_t}$ by \eq{PH}.  The last part follows by Definition \ref{s6dfn5}.
\end{proof}

\begin{notes} Recall the notation of Definitions \ref{s6dfn1} and \ref{s6dfn5} and Propositions \ref{s6prop1} and \ref{s6prop2}.
\begin{itemize} 
\item[(a)] Since $dH_t=F_t$, $d\psi(t)=\varphi|_{\tilde{N}(t)}$ and $d\alpha_{H_t}=d\alpha$, we see that $dP_{H_t}=P_{F_t}$. 
\item[(b)] The notation introduced for the maps $F_t$, $P_{F_t}$, $H_t$ and $P_{H_t}$ and for the form $\psi(t)$ will be
used throughout the remaining sections.
\end{itemize}
\end{notes}

\section{Sobolev embedding inequalities and\\ regularity}\label{sei}

The most important analytic results of this paper are contained in this section.  The first set are known as \emph{Sobolev embedding
inequalities}.  These are estimates for the norm of a form in terms of the norm of its exterior derivative, both in suitable Sobolev spaces.  
The most challenging problem is to obtain an estimate on $\tilde{N}(t)$ which is \emph{independent} of $t$, which is where we shall begin.  Later,
we look at \emph{elliptic regularity} and see how it may be applied to show that solutions to $F_t(\alpha)=0$ are smooth, given that $\alpha$ has 
enough initial differentiability.

\subsection{Sobolev embedding inequalities on {\boldmath $\tilde{N}(t)$}}\label{seis1}

Recall the notation of $\S$\ref{s4subs1}-\ref{s4subs2}.  
We begin this subsection with the derivation of Sobolev embedding inequalities on $A$ and $\hat{N}$ which are for `\emph{unweighted}' 
Sobolev spaces -- this
is essential for our application.
Moreover, so that these inequalities are useful for our purposes, 
we first need to construct finite-dimensional spaces of compactly supported forms that `approximate' certain self-dual
harmonic forms on $A$ and $\hat{N}$.

\begin{prop}\label{seiprop1}  Recall the notation of Definition \ref{ch6s2dfn1}
and $\S$\ref{s4subs1}(c) and let
$g_0$ be the Euclidean metric on $\R^7$.
There is a finite-dimensional subspace $\mathcal{K}_{\text{\emph{ap}}}^{A}$ of\/ $C^1_{\text{\emph{cs}}}\big((\Lambda^2_+)_{g_0|_A}T^*A\big)$ 
which is isomorphic to the kernel $\mathcal{K}^{A}$ of 
$$d:L^{4/3}_{1,\,-2}\big((\Lambda^2_+)_{g_0|_A}T^*A\big)\rightarrow
 L^{4/3}_{0,\,-3}(\Lambda^3T^*A)$$
such that if $\alpha\in\mathcal{K}^{A}$ satisfies $\langle\alpha,\beta\rangle_{L^2}=0$ for all
$\beta\in\mathcal{K}_{\text{\emph{ap}}}^{A}$, $\alpha=0$.
\end{prop}

\begin{proof}  Let $\rho_A$ be a radius function on $A$, as given in Definition \ref{ch6s1dfn2}.
The map
\begin{equation}\label{sei1}d+d^*:L^{4/3}_{1,\,-2}\big((\Lambda^2_+)_{g_0|_A}T^*A\oplus\Lambda^4T^*A\big)\rightarrow 
L^{4/3}_{0,\,-3}(\Lambda^3T^*A)
\end{equation}
is elliptic and therefore has finite-dimensional kernel.  Moreover, if $(\alpha,\beta)$ lies in the kernel, we see that $dd^*\beta=\Delta\beta=0$,
and so $*\beta$ is a harmonic function which is certainly of order $O(\rho_A^{-2})$ as $\rho_A\rightarrow\infty$.  The maximum principle for harmonic
functions allows us to deduce that $*\beta$, hence $\beta$, is zero.  Thus $\mathcal{K}^{A}$ is isomorphic to the kernel of \eq{sei1} and
is finite-dimensional.  Since $C^1_{\text{cs}}$ is dense in $L^{4/3}_{1,\,-2}$ we may choose $\mathcal{K}^{A}_{\text{ap}}$ 
as claimed.
\end{proof}

\begin{prop}\label{seiprop2} Use the notation of Definition \ref{ch6s2dfn1}, Lemma \ref{s3sub4lem1}, $\S$\ref{s4subs1}(c) and Proposition \ref{seiprop1}.
Let $\alpha\in L^{4/3}_{1,\,-2}\big((\Lambda^2_+)_{g_0|_A}T^*A\big)$ be such that $\langle\alpha,\beta\rangle_{L^2}=0$ for all 
$\beta\in\mathcal{K}_{\text{\emph{ap}}}^{A}$.  There exists $C(A)>0$ independent of $\alpha$ such that
$$\|\alpha\|_{L^2}\leq C(A)\|d\alpha\|_{L^{4/3}}.$$
\end{prop}
  
\begin{proof}
First note that $L^{4/3}_{1,\,-2}$ embeds continuously in $L^2_{0,\,-2}=L^2$ by Theorem \ref{ch6s2thm1}, 
using the observation made after Definition
\ref{ch6s2dfn1}. Thus, there exists a constant $C(A)^{\prime}>0$, independent of $\alpha$, such that
$$\|\alpha\|_{L^2}\leq C(A)^{\prime}\|\alpha\|_{L^{4/3}_{1,\,-2}}.$$ 
Let $\beta\in L^{4/3}_{1,\,-2}(\Lambda^4T^*A)$. 
By Proposition \ref{s3sub4prop3} and Lemma \ref{s3sub4lem1},
if $\alpha$ is $L^2$-orthogonal to $\mathcal{K}^{A}$, 
 there exists a constant $C(A)^{\prime\prime}>0$, independent of $(\alpha,\beta)$, such that
$$\|(\alpha,\beta)\|_{L^{4/3}_{1,\,-2}}\leq C(A)^{\prime\prime}\|d\alpha+d^*\beta\|_{L^{4/3}_{0,\,-3}}=C(A)^{\prime\prime}\|d\alpha+d^*\beta\|_{L^{4/3}},$$
again using the note after Definition \ref{ch6s2dfn1}.
However, by the choice of $\mathcal{K}^A_{\text{ap}}$, this inequality still holds for $\alpha$ $L^2$-orthogonal to $\mathcal{K}^A_{\text{ap}}$, 
increasing $C(A)^{\prime\prime}$ if necessary.
Setting $\beta=0$ and using our earlier relationship between norms of $\alpha$, we get our result with $C(A)=C(A)^{\prime}C(A)^{\prime\prime}$.
\end{proof}

\begin{note}
The Sobolev embedding inequality given in the proposition above is invariant under dilations of $A$.  This is key for our purposes.
\end{note}

Having constructed an `approximation' to the self-dual harmonic forms on $A$ in $L^{4/3}_{1,\,-2}$ and proved a Sobolev embedding inequality
related to it, we now give analogous results for $\hat{N}$.

\begin{prop}\label{seiprop3} Recall the notation in Definition \ref{ch6s2dfn1} and $\S$\ref{s4subs1}(b).
There is a finite-dimensional subspace $\mathcal{K}_{\text{\emph{ap}}}^{N}$ of\/
 $C^1_{\text{\emph{cs}}}\big((\Lambda^2_+)_{g|_{\hat{N}}}T^*\hat{N}\big)$ 
which is isomorphic to the kernel $\mathcal{K}^{N}$ of 
$$d:L^{4/3}_{1,\,-2}\big((\Lambda^2_+)_{g|_{\hat{N}}}T^*\hat{N}\big)\rightarrow L^{4/3}_{0,\,-3}(\Lambda^3T^*\hat{N}),$$
such that if $\alpha\in\mathcal{K}^{N}$ satisfies $\langle\alpha,\beta\rangle_{L^2}=0$ for all
$\beta\in\mathcal{K}_{\text{\emph{ap}}}^{N}$, $\alpha=0$.
\end{prop}

\begin{proof}
The one point in the proof of Proposition \ref{seiprop1} that cannot be mirrored is the use of the maximum
principle argument.  Instead note that,
if $(\alpha,\beta)$ lies in the kernel of \eq{sei1}, now acting on $\hat{N}$, it is clear
through an integration by parts argument that $d\alpha=d^*\beta=0$, valid since 
$L^{4/3}_{1,\,-2}\hookrightarrow L^2_{0,\,-2}$ by
Theorem \ref{ch6s2thm1}.  Therefore, $\mathcal{K}^N$ can be considered as a subspace of the kernel of \eq{sei1} and so is finite-dimensional.
\end{proof}

\noindent The next result is proved in an entirely similar manner to Proposition \ref{seiprop2}.

\begin{prop}\label{seiprop4} Use the notation of Definition \ref{ch6s2dfn1}, $\S$\ref{s4subs1}(b) and Proposition \ref{seiprop3}.
Let $\alpha\in L^{4/3}_{1,\,-2}\big((\Lambda^2_+)_{g|_{\hat{N}}}T^*\hat{N}\big)$ be such that $\langle\alpha,\beta\rangle_{L^2}=0$ for all 
$\beta\in\mathcal{K}_{\text{\emph{ap}}}^{N}$.  There exists $C(N)>0$ independent of $\alpha$ such that
$$\|\alpha\|_{L^2}\leq C(N)\|d\alpha\|_{L^{4/3}}.$$
\end{prop}

So that we can derive our main Sobolev embedding inequality on $\tilde{N}(t)$,
we construct a suitable `approximation' to the $L^2$ self-dual harmonic forms on $\tilde{N}(t)$ 
from $\mathcal{K}^A_{\text{ap}}$ and $\mathcal{K}^N_{\text{ap}}$.  

\begin{dfn}\label{seidfn1}  Use the notation of $\S$\ref{s4subs1}(a)-(c), Definition \ref{s4dfn2} and
Propositions \ref{seiprop1} and \ref{seiprop3}.  In particular, recall we fixed $\nu\in(0,1)$.

Define a diffeomorphism $\tilde{\Psi}_A$ between
 $tY_A(t)$ and $\chi(tK_A)\sqcup\Psi_{\tilde{N}(t)}\big((tR,\epsilon)\times\Sigma\big)$ by $\tilde{\Psi}_A(x)=\chi(x)$ for all $x\in tK_A$ and
 $\tilde{\Psi}_A\big(t\Phi_A(r,\sigma)\big)=\Psi_{\tilde{N}(t)}(tr,\sigma)$ for all $(r,\sigma)\in(R,t^{-1}\epsilon)\times\Sigma$.  
If $\tau$ is sufficiently small
 we may identify the metrics, hence the self-dual 2-forms, on $tY_A(t)$ and $\tilde{\Psi}_A\big(tY_A(t)\big)$.

Let $\alpha_1,\ldots,\alpha_l$ be a basis for $\mathcal{K}^A_{\text{ap}}$.
Since $t^{\nu-1}\rightarrow\infty$ as $t\rightarrow 0$, 
by choosing $\tau$ sufficiently small we may ensure that
 the support of $\alpha_j$ is contained in $K_A\sqcup\Phi_A\big((R,t^{\nu-1})\times\Sigma\big)\subseteq Y_A(t)$ for $j=1,\ldots,l$.
Denote by $t\alpha_j$ the corresponding form to $\alpha_j$ on the dilated submanifold $tY_A(t)$.  Let 
$\xi_j\in C^1\big((\Lambda^2_+)_{\tilde{g}(t)}T^*\tilde{N}(t)\big)$
be equal to zero outside of $\tilde{\Psi}_A\big(tK_A\sqcup t\Phi_A\big((R,t^{\nu-1})\times\Sigma\big)\big)$ 
and such that $\tilde{\Psi}_A^*(\xi_j)=t\alpha_j$, for $j=1,\ldots,l$.
This is possible by the identifications made above for sufficiently small $\tau$.

Let $0<\nu^{\prime\prime}<\nu^{\prime}<\nu$.  Then $t^{\nu}<t^{\nu^{\prime}}<t^{\nu^{\prime\prime}}<\epsilon$ for all $t\in(0,\tau)$,
if $\tau$ is small enough.  Let $\xi_{l+1},\ldots,\xi_m$ be a basis for $\mathcal{K}^N_{\text{ap}}$.  
Since $t^{\nu^{\prime\prime}}\rightarrow0$ as $t\rightarrow 0$, it is possible to be certain that,
 if $\tau$ is made smaller as necessary, the support of $\xi_j$ is contained in 
$K_N\sqcup\Psi_N\big((t^{\nu^{\prime\prime}},\epsilon)\times\Sigma\big)$ for $j=l+1,\ldots,m$.  
We may thus extend by zero each $\xi_j$ to a $C^1$ 
self-dual 2-form on $\tilde{N}(t)$, with respect to $\tilde{g}(t)$, since $\tilde{g}(t)$ is conformal to $g|_{\hat{N}}$ over
 $K_N\sqcup\Psi_N\big((t^{\nu},\epsilon)\times\Sigma\big)$.

Let $\tilde{\mathcal{K}}_{\text{ap}}(t)$ be the subspace of $C^1\big((\Lambda^2_+)_{\tilde{g}(t)}T^*\tilde{N}(t)\big)$ 
spanned by $\xi_1,\ldots,\xi_m$.  
We have thus constructed a space of self-dual 2-forms on $\tilde{N}(t)$ 
such that the spaces $\text{Span}\{\xi_1,\ldots,\xi_l\}$ and $\text{Span}\{\xi_{l+1},\ldots,\xi_m\}$
are naturally identified with $\mathcal{K}^A_{\text{ap}}$ and $\mathcal{K}^N_{\text{ap}}$ respectively.
\end{dfn}

\begin{remark} 
The notion that $\tilde{\mathcal{K}}_{\text{ap}}(t)$ `approximates' the self-dual harmonic forms on $\tilde{N}(t)$ in $L^2$
 will be explored later, specifically in Propositions \ref{seiprop6} and Corollary \ref{seicor2}, but it 
highlights the similarity between the main results in $\S$\ref{seis1}-\ref{seis2} and \cite[Proposition 4.2]{Kovalev2}.  
\end{remark}

\begin{thm}\label{seithm1} Recall the notation of Definitions \ref{s4dfn2} and \ref{seidfn1}.
 Let $\alpha$ be an element of $L^{4/3}_1\big((\Lambda^2_+)_{\tilde{g}(t)}T^*\tilde{N}(t)\big)$ such that $\langle\alpha,\beta\rangle_{L^2}=0$ for all
$\beta\in\tilde{\mathcal{K}}_{\text{\emph{ap}}}(t)$. 
There is a constant $C(\tilde{N})>0$, independent of $\alpha$ and $t$, such that
$$\|\alpha\|_{L^2}\leq C(\tilde{N})\|d\alpha\|_{L^{4/3}}.$$
\end{thm}

\begin{proof}  We shall use the notation of Propositions \ref{seiprop1}-\ref{seiprop4} and Definitions \ref{ch6s2dfn1}, 
\ref{s4dfn2} and \ref{seidfn1} in this proof.
  
Let
$f_{\text{dec}}:\R\rightarrow [0,1]$ be a smooth decreasing function such that
$$f_{\text{dec}}(x)=\left\{\begin{array}{lll} 1 & & \text{for $x\leq \nu^{\prime\prime}$,}\\
0&&\text{for $x\geq\nu^{\prime}$}\end{array}\right.$$ and
$f_{\text{dec}}(x)\in(0,1)$ for $x\in(\nu^{\prime\prime},\nu^{\prime})$.  Define $f:\tilde{N}(t)\rightarrow [0,1]$ to be the smooth function
given by
$$f(x)=\left\{\begin{array}{lll}  0 & & \text{for $x\in\chi(tK_A)$,}\\[2pt]
f_{\text{dec}}\left(\frac{\log r}{\log t}\right) & & \text{for $x=\Psi_{\tilde{N}(t)}(r,\sigma)
\in\Psi_{\tilde{N}(t)}\big((tR,\epsilon)\times\Sigma\big)$ and}\\[4pt]
1 & & \text{for $x\in K_N$.}\end{array}\right.$$
Consider $\alpha=f\alpha+(1-f)\alpha$.

Since the support of $f$ is contained in $\tilde{N}_u(t)$, 
which is a subset of $\hat{N}$, $f\alpha$ can be extended by zero to the rest
of $\hat{N}$ to define a compactly supported element of $L^{4/3}_{1,\,-2}\big((\Lambda^2_+)_{g|_{\hat{N}}}T^*\hat{N}\big)$.  
If $\beta\in\mathcal{K}^N_{\text{ap}}$, it may 
be written as a linear combination of $\xi_{l+1},\ldots,\xi_m\in\tilde{\mathcal{K}}_{\text{ap}}(t)$ and thus, as the support of $\beta$ is contained 
in a subset of $\tilde{N}(t)$ upon which $f=1$,
$$\langle f\alpha,\beta \rangle_{L^2}=\langle \alpha , \beta \rangle_{L^2}=0.$$
Applying Proposition \ref{seiprop4}, noting that the metrics $\tilde{g}(t)$ and $g|_{\hat{N}}$ are conformal over $\tilde{N}_u(t)$, gives that
\begin{align*}
\|f\alpha\|_{L^2}&\leq C(N) \|d(f\alpha)\|_{L^{4/3}}\\
&\leq C(N)\big(\|fd\alpha\|_{L^{4/3}}+\|df\w\alpha\|_{L^{4/3}}\big)\\
&\leq C(N)\|d\alpha\|_{L^{4/3}}+C(N)\|df\|_{L^4}\|\alpha\|_{L^2},
\end{align*}
using H\"older's inequality and the fact that $f(x)\in [0,1]$ for all $x\in\tilde{N}(t)$.

In a similar vein, the support of $(1-f)\alpha$ is contained in $\tilde{\Psi}_A\big(tY_A(t)\big)$.  It may then be identified, using $\tilde{\Psi}_A$, 
with a compactly supported element, $\gamma$ say,
of $L^{4/3}_{1,\,-2}\big((\Lambda^2_+)_{g_0|_{tA}}T^*(tA)\big)$, where $g_0$ is the Euclidean metric on $\R^7$.  
Moreover, if $\beta\in\mathcal{K}^A_{\text{ap}}$, it may be
written as $\beta=\sum_{j=1}^la_j\alpha_j$, where the $\alpha_j$ are given in Definition \ref{seidfn1}.  
Thus, the corresponding form $t\beta$ on $tY_A(t)$ is identified with $\xi=\sum_{j=1}^la_j\xi_j$ as
described in 
Definition \ref{seidfn1}.  Therefore, as the support of $\xi$ is contained in the subset of $\tilde{N}(t)$ upon which $f=0$,
and we identify the metrics on $tY_A(t)$ and $\tilde{\Psi}_A\big(tY_A(t)\big)$,
$$\langle\gamma,t\beta\rangle_{L^2}=\langle(1-f)\alpha,\xi\rangle_{L^2}=\langle\alpha,\xi\rangle_{L^2}=0.$$
We may thus apply Proposition \ref{seiprop2} to $\gamma$, recalling that the inequality is invariant under dilations, and see that
$$\|\gamma\|_{L^2}\leq C(A)\|d\gamma\|_{L^{4/3}},$$
where the norms are calculated using $g_0|_{tA}$.  By increasing $C(A)$ to $2C(A)$ we can ensure, for small $\tau$, that the inequality holds with
norms calculated using $\tilde{g}(t)$.  We may therefore deduce that
$$\|(1-f)\alpha\|_{L^2}\leq 2C(A)\|d((1-f)\alpha)\|_{L^{4/3}}\leq 2C(A)\big(\|d\alpha\|_{L^{4/3}}+\|d(1-f)\|_{L^4}\|\alpha\|_{L^2}\big),$$
using the same method as earlier.

We must now estimate $\|df\|_{L^4}=\|d(1-f)\|_{L^4}$.  It is clear that $df$ is zero except on 
$\Psi_{\tilde{N}(t)}\big((t^{\nu^{\prime}},t^{\nu^{\prime\prime}})\times\Sigma\big)$, where
$$\Psi_{\tilde{N}(t)}^*(df)=\frac{1}{r\log t}\frac{df}{dr}\left(\frac{\log r}{\log t}\right)dr.$$
Hence, the $L^4$ norm of $df$ is of the same order as
$$\left(\int^{t^{\nu^{\prime\prime}}}_{t^{\nu^{\prime}}} \frac{1}{(r\log t)^4}r^3dr\right)^{\frac{1}{4}}=
(\nu^{\prime}-\nu^{\prime\prime})^{\frac{1}{4}}|\log t|^{-\frac{3}{4}},$$
which tends to zero as $t\rightarrow 0$.  

Combining our estimates we see that
$$\|\alpha\|_{L^2}\leq\|f\alpha\|_{L^{2}}+\|(1-f)\alpha\|_{L^{2}}\leq \big(C(N)+2C(A)\big)\big(\|d\alpha\|_{L^{4/3}}+\|df\|_{L^4}\|\alpha\|_{L^2}\big)$$
and, by rearranging, that
$$\big(1-\big(C(N)+2C(A)\big)\|df\|_{L^4}\big)\|\alpha\|_{L^2}\leq \big(C(N)+2C(A)\big)\|d\alpha\|_{L^{4/3}}.$$
Since $\|df\|_{L^4}\rightarrow 0$ as $t\rightarrow 0$, we may choose $\tau$ sufficiently small that 
$$\big(1-(C(N)+2C(A)\big)\|df\|_{L^4})\geq\frac{1}{2}$$
for all $t\in (0,\tau)$.  Having made this choice, our result holds with $C(\tilde{N})=2\big(C(N)+2C(A)\big)$.
\end{proof}

\noindent This theorem is an analogue of \cite[Theorem 6.12]{Joyce3} and the proof presented is essentially an adaptation of the 
proof given in the reference cited.

\subsection{Further inequalities and observations}\label{seis2}

We build upon the Sobolev embedding inequality given in Theorem \ref{seithm1} so that we can estimate the
 $L^{4/3}_{1,\,-2,\,t}$ and $L^8_{2,\,1,\,t}$ 
norms of a self-dual 2-form $\alpha$ by the $L^{4/3}$ and $L^8_{1,\,0,\,t}$ norms of $d\alpha$.  These inequalities will depend upon $t$, so we have to understand 
how it varies with $t$, particularly as $t\rightarrow 0$.

\begin{lemma}\label{seilem1}  Use the notation of Definitions \ref{s4dfn2} and \ref{s4dfn3}.
If $\tau$ is sufficiently small, there exists a constant $C(\rho)>0$, independent of $t$, such that the radius function $\rho_t$ on $\tilde{N}(t)$ satisfies
\begin{align*}
\|\rho_t^{-2}\|_{L^2}\leq C(\rho)|\log t|^{\frac{1}{2}}
\quad\text{and}\quad\|\rho_t^{-5}\|_{L^2}\leq C(\rho)t^{-3}.
\end{align*}
\end{lemma}

\begin{proof}
This is a straightforward calculation.  Recall that $\tilde{g}(t)$ is the chosen metric on $\tilde{N}(t)$ as in Definition \ref{s4dfn4}.

Consider $\|\rho_t^{-2}\|_{L^2}$.
  On $K_N$, $\rho_t=1$ and so 
$$\int_{K_N}\rho_t^{-4}dV_{\tilde{g}(t)}=\text{vol}(K_N)=O(1).$$  On $\chi(tK_A)$, $\rho_t=t$ and $\text{vol}(tK_A)=O(t^4)$, so
$$\int_{\chi(tK_A)}\rho_t^{-4}dV_{\tilde{g}(t)}=O(t^{-4}\text{vol}(tK_A))=O(1)$$
as well.  Our main contribution to $\|\rho_t^{-2}\|_{L^2}$ therefore comes from the integral
 over $\Psi_{\tilde{N}(t)}\big((tR,\epsilon)\times\Sigma\big)$.
By definition, $\Psi_{\tilde{N}(t)}^*(\rho_t^{-2})=O(r^{-2})$ and
$$\int_{(tR,\epsilon)\times\Sigma}r^{-4}dV_{g_{\text{cone}}}=O\left(\int_{tR}^{\epsilon}r^{-1}dr\right)=O(|\log t|).$$
Hence, the first inequality follows.

Now consider $\|\rho_t^{-5}\|_{L^2}$.  Again the integral over $K_N$ is $O(1)$ but  
$$\int_{\chi(tK_A)}\rho_t^{-10}dV_{\tilde{g}(t)}=O\big(t^{-10}\text{vol}(tK_A)\big)=O(t^{-6}).$$
The integral over $\Psi_{\tilde{N}(t)}\big((tR,\epsilon)\times\Sigma\big)$ is of the same order as
$$\int_{(tR,\epsilon)\times\Sigma}r^{-10}dV_{g_{\text{cone}}}=O\left(\int_{tR}^{\epsilon}r^{-7}dr\right)=O(t^{-6}).$$
From these calculations we deduce the second inequality.
\end{proof}

\begin{prop}\label{seilem2} Use the notation of Definitions \ref{s4dfn2}-\ref{s4dfn4} and \ref{seidfn1}.
Let $\alpha\in L^{8}_{2}\big((\Lambda^2_+)_{\tilde{g}(t)}T^*\tilde{N}(t)\big)$ 
be such that $\langle\alpha,\beta\rangle_{L^2}=0$ for all 
$\beta\in\tilde{\mathcal{K}}_{\text{\emph{ap}}}(t)$. 
If $\tau$ is sufficiently small, there exists a constant $C(\tilde{N},\rho)>0$, independent of $\alpha$ and $t$, such that
\begin{align*}
\|\alpha\|_{L^{4/3}_{1,\,-2,\,t}}&\leq C(\tilde{N},\rho)|\log t|^{\frac{1}{2}}\|d\alpha\|_{L^{4/3}}\\
\intertext{and}
\|\alpha\|_{L^8_{2,\,1,\,t}}&\leq C(\tilde{N},\rho)\left(\|d\alpha\|_{L^8_{1,\,0,\,t}}+t^{-3}\|d\alpha\|_{L^{4/3}}\right)
\end{align*}
\end{prop}

\begin{proof}
By Proposition \ref{s5prop8} with $p=4/3$, $k=0$ and $\eta=-2$, 
there exists a constant $C(d+d^*)>0$, independent of $\alpha$ and $t$, such that
$$\|\alpha\|_{L^{4/3}_{1,\,-2,\,t}}\leq C(d+d^*)\left(\|d\alpha\|_{L^{4/3}_{0,\,-3,\,t}}
+\|\alpha\|_{L^1_{0,\,-2,\,t}}\right).$$
Note that
$$\|\alpha\|_{L^1_{0,\,-2,\,t}}=\|\rho_t^{-2}\alpha\|_{L^1}\leq\|\rho_t^{-2}\|_{L^2}\|\alpha\|_{L^2}$$
by H\"older's inequality.  Applying Theorem \ref{seithm1} and Lemma \ref{seilem1} we see that
$$\|\alpha\|_{L^1_{0,\,-2,\,t}}\leq C(\rho)C(\tilde{N})|\log t|^{\frac{1}{2}}\|d\alpha\|_{L^{4/3}}.$$
Since the $L^{4/3}_{0,\,-3,\,t}$ norm equals the $L^{4/3}$ norm, as can be seen by Definition \ref{s4dfn4}, the first result follows.

We can similarly appeal to Proposition \ref{s5prop8}, Theorem \ref{seithm1} and Lemma \ref{seilem1} to deduce the second inequality.
\end{proof}

We now need to understand the relationship better between $\tilde{\mathcal{K}}_{\text{ap}}(t)$ and the harmonic self-dual 2-forms on $\tilde{N}(t)$.

\begin{dfn}\label{H2+} Recall the notation of Definitions \ref{s4dfn2}-\ref{s4dfn4}.
Let 
$$\tilde{\mathcal{H}}^2(t)=\big\{\alpha\in L^2\big(\Lambda^2T^*\tilde{N}(t)\big)\,:\,d\alpha=d^*\alpha=0\big\}$$
and let
$$\tilde{\mathcal{H}}^2_+(t)=\tilde{\mathcal{H}}^2(t)\cap C^{\infty}\big((\Lambda^2_+)_{\tilde{g}(t)}T^*\tilde{N}(t)\big).$$
Note that $\text{dim}\,\tilde{\mathcal{H}}^2_+(t)=b^2_+\big(\tilde{N}(t)\big)$ by Hodge theory.
\end{dfn}

\noindent 
We show that the $L^2$-orthogonal projection from $\tilde{\mathcal{H}}^2_+(t)$ to $\tilde{\mathcal{K}}_{\text{ap}}(t)$ is
injective.

\begin{prop}\label{seiprop6}
If $\alpha\in\tilde{\mathcal{H}}^2_+(t)$ satisfies $\langle\alpha,\beta\rangle_{L^2}=0$ for all $\beta\in\tilde{\mathcal{K}}_{\text{\emph{ap}}}(t)$,
$\alpha=0$, using the notation of Definitions \ref{seidfn1} and \ref{H2+}.
\end{prop}

\begin{proof}
If $\alpha\in\tilde{\mathcal{H}}^2_+(t)$, $\alpha\in L^{4/3}_1\big((\Lambda^2_+)_{\tilde{g}(t)}T^*\tilde{N}(t)\big)$
since $\alpha$ is smooth and $\tilde{N}(t)$ is compact.   
If, additionally, $\langle\alpha,\beta\rangle_{L^2}=0$, we may apply Theorem \ref{seithm1}
to deduce that $\|\alpha\|_{L^2}=0$, because $d\alpha=0$, and thus $\alpha=0$.
\end{proof}

\noindent We can now state a useful result that follows immediately from Corollary \ref{seqcor} and Proposition \ref{seiprop6}.

\begin{cor}\label{seicor2}
The space $\tilde{\mathcal{H}}^2_+(t)$ is isomorphic to $\tilde{\mathcal{K}}_{\text{\emph{ap}}}(t)$ via $L^2$-orthogonal projection,
in the notation of Definitions \ref{seidfn1} and \ref{H2+}.
\end{cor}

%

Our final result in this subsection is a useful refinement of Corollary \ref{s2cor1}.

\begin{cor}\label{seicor1}
Use the notation of Definitions \ref{s4dfn2}-\ref{s4dfn4} and \ref{seidfn1}.  Let 
\begin{align*}
V^8_{2}(t)=\big\{\alpha\in L^8_{2}\big((\Lambda^2_+)_{\tilde{g}(t)}T^*\tilde{N}(t)\big)\,:\,\langle\alpha,\beta\rangle_{L^2}=0\;\text{for all 
 $\beta\in\tilde{\mathcal{K}}_{\text{\emph{ap}}}(t)$}\big\}
\end{align*}
and let $W^8_1(t)$ denote the space of $L^8_1$ exact 3-forms on $\tilde{N}(t)$.
The map $d:V^8_{2}(t)\rightarrow W^8_1(t)$ is a linear isomorphism between Banach spaces, which therefore has
a linear inverse $d^{-1}:W^8_1(t)\rightarrow V^8_2(t)$.
\end{cor}

\begin{proof}
For injectivity, let $\alpha\in V^8_{2}(t)$ satisfy $d\alpha=0$.
 Then $\alpha\in\tilde{\mathcal{H}}^2_+(t)$, as given in Definition \ref{H2+}.  By Proposition \ref{seiprop6}, $\alpha=0$.
Since the exterior derivative is clearly linear and $W^8_1(t)$ is shown to be a Banach space in the proof of Proposition \ref{s2prop3}, 
we need only check that $V^8_{2}(t)$ is a Banach space and that 
the map is surjective.  The first point is clear since the condition on $\alpha\in L^8_{2}$ is a closed one.
For surjectivity, let $\gamma\in W^8_1(t)$.
By Proposition \ref{s2prop3}, 
there exists $\xi\in L^8_{2}\big((\Lambda^2_+)_{\tilde{g}(t)}T^*\tilde{N}(t)\big)$ such that $d\xi=\gamma$. 
Now, we may certainly take an $L^2$-orthogonal projection of $\xi$ and, by Corollary \ref{seicor2},
 get a form $\alpha\in V^8_{2}(t)$ such that $d\alpha=d\xi=\gamma$.
\end{proof}

\begin{dfn}\label{seidfn3} Use the notation of Proposition \ref{s5prop3} and Corollary \ref{seicor1} and recall the maps $F_t$ and $H_t$
given in Definitions \ref{s6dfn1} and \ref{s6dfn5}.
We shall abuse notation and write $H_t(\alpha)=d^{-1}\big(F_t(\alpha)\big)\in V^8_{2}(t)$ 
for $\alpha\in L^8_{2}\big(\tilde{U}(t)\big)$
with $\|\alpha\|_{C^1_{1,\,t}}$ sufficiently small.  We shall adopt similar conventions
for $\psi(t)$, $\alpha_{H_t}$ and $P_{H_t}$, as previously given in Definition \ref{s6dfn5}.  
The properties of $H_t$, $\psi(t)$ and $P_{H_t}$ described thus far
still hold under this change.  
\end{dfn}

\subsection{Elliptic regularity}

Recall the deformation map $F_t$ and the open neighbourhood $\tilde{U}(t)$ of the zero section in $(\Lambda^2_+)_{\tilde{g}(t)}T^*\tilde{N}(t)$.
Let $\alpha\in L^p_{k+1}\big(\tilde{U}(t)\big)$ for some choice of
$p> 4$ and $k\geq 1$.
The key to showing that solutions of $F_t(\alpha)=0$ are smooth  
is to use \emph{elliptic regularity} results.

\begin{dfn}\label{s8dfn}  Use the notation of Proposition \ref{s5prop3} and Definition \ref{s6dfn1}.  
Let $p>4$, let $k\geq 1$, let $l\in\N$ and let $a\in(0,1)$.
For $\alpha\in L^p_{k+1}\big(\tilde{U}(t)\big)$  or $C^{l+1,\,a}\big(\tilde{U}(t)\big)$, define 
$$G_t(\alpha)=\pi_{\Lambda^2_+}\big(d^*F_t(\alpha)\big),$$ 
where $\pi_{\Lambda^2_+}$ is the projection from 2-forms to self-dual 2-forms.  Since $\alpha$ is self-dual,
$G_t(\alpha)$ is elliptic in $\alpha$ at zero, meaning that its linearisation at zero is elliptic.  
\end{dfn}

We can now describe the map $G_t$ further.

\begin{dfn}\label{s8dfn2}  Use the notation of Propositions \ref{s5prop3} and Definitions 
\ref{s4dfn2}-\ref{s4dfn4}, \ref{s6dfn1} and \ref{s8dfn}.  
Let $\alpha\in L^{p}_{k+1}\big(\tilde{U}(t)\big)$ or $C^{l+1,\,a}\big(\tilde{U}(t)\big)$.  We can write, for $x\in\tilde{N}(t)$,
$$G_t(\alpha)(x)=R_{F_t}\big(x,\alpha(x),\nabla_t\alpha(x)\big)\nabla_t^2\alpha(x)+E_{F_t}\big(x,\alpha(x),\nabla_t\alpha(x)\big),$$
where $R_{F_t}$ and $E_{F_t}$ are smooth functions of their arguments.  
This is possible since $G_t$ is linear in $\nabla_t^2\alpha$ with coefficients depending on $\alpha$ and $\nabla_t\alpha$.  Let $R_{F_t}(\alpha)(x)=R_{F_t}\big(x,\alpha(x),\nabla_t\alpha(x)\big)$ 
and make a similar definition for
$E_{F_t}(\alpha)$.  We see that the \emph{linear} map defined for 
self-dual 2-forms $\gamma$ on $\tilde{N}(t)$ by
$$S_{\alpha,\,t}(\gamma)=R_{F_t}(\alpha)\nabla_t^2\gamma$$
is elliptic at zero.  Note that $S_{\alpha,\,t}$ is \emph{not} the linearisation.  
Moreover, by Proposition \ref{s6prop1}, 
$E_{F_t}(\alpha)$ is the sum of a constant term, $\pi_{\Lambda^2_+}\big(d^*\varphi|_{\tilde{N}(t)}\big)$, and a term which is 
at worst quadratic in $\rho_t^{-1}\alpha$ and $\nabla_t\alpha$.  Thus 
$E_{F_t}(\alpha)\in L^{p}_k\big((\Lambda^2_+)_{\tilde{g}(t)}T^*\tilde{N}(t)\big)$ or $C^{l,\,a}\big((\Lambda^2_+)_{\tilde{g}(t)}T^*\tilde{N}(t)\big)$,
if $\alpha$ is small enough in the $C^1_{1,\,t}$ norm. 
\end{dfn}

We now employ the standard ``bootstrap'' method to prove the regularity result we require.  

\begin{prop}\label{s8prop2} Use the notation of Proposition \ref{s5prop3} and Definitions \ref{s4dfn4} and \ref{s6dfn1}.
If $\alpha\in L^p_2\big(\tilde{U}(t)\big)$, for $p>4$, satisfies $F_t(\alpha)=0$ and $\|\alpha\|_{C^1_{1,\,t}}$ is sufficiently small, 
$\alpha\in C^{\infty}\big(\tilde{U}(t)\big)$.
\end{prop}

\begin{proof}
Recall the notation and observations in Definition \ref{s8dfn2}.
Suppose $\alpha\in C^{k+1,\,a}\big(\tilde{U}(t)\big)$ for some $k\in\N$ and $a=1-4/p\in(0,1)$.  
Suppose also that
 $F_t(\alpha)=0$,
$\|\alpha\|_{C^1_{1,\,t}}$ is sufficiently small, $\gamma\in C^{k+1,\,a}\big((\Lambda^2_+)_{\tilde{g}(t)}T^*\tilde{N}(t)\big)$ and 
$S_{\alpha,\,t}(\gamma)\in C^{k,\,a}$.  
The Schauder estimates, which relate to linear elliptic operators between H\"older spaces,
 are then applicable since $S_{\alpha,\,t}$ is a \emph{linear} second-order elliptic operator whose coefficients lie in $C^{k,\,a}$.
We may thus deduce that $\gamma\in C^{k+2,\,a}\big((\Lambda^2_+)_{\tilde{g}(t)}T^*\tilde{N}(t)\big)$.  Notice that
$$S_{\alpha,\,t}(\alpha)=-E_{F_t}(\alpha)\in C^{k,\,a}\big((\Lambda^2_+)_{\tilde{g}(t)}T^*\tilde{N}(t)\big),$$
as noted in Definition \ref{s8dfn2}.  Hence, 
we can set $\gamma=\alpha$ and see that $\alpha$ lies in $C^{k+2,\,a}$, knowing only a priori that $\alpha$ lay in $C^{k+1,\,a}$.  

Since $L^p_2$ embeds in $C^{1,\,a}$ by the Sobolev Embedding Theorem, the result follows by induction.
\end{proof}

\section{Desingularization: stage 2}\label{s7}

Here we tackle the bulk of the analysis, using the foundational work of the last two sections, which will enable us to prove our main result 
in $\S$\ref{s9}.  We first estimate the size of $\varphi$ on various parts of $\tilde{N}(t)$, then use this estimate to bound the
modulus of $\psi(t)$.  We briefly derive some estimates on $P_{H_t}$ before constructing a \emph{contraction map}.  The fixed point $\alpha(t)$ of this
contraction satisfies $F_t\big(\alpha(t)\big)=0$ and so is a candidate for a form corresponding to a coassociative deformation of $\tilde{N}(t)$.  

\subsection{Estimating {\boldmath $\varphi|_{\tilde{N}(t)}$}}

\begin{prop}\label{s7prop1} Recall the notation introduced in $\S$\ref{s4subs1}(a)-(c) and Definition \ref{s4dfn2}.  
In particular, recall that we have a constant $\nu\in(0,1)$, a diffeomorphism $\chi:B(0;\epsilon_M)\rightarrow V\subseteq M$ and a
diffeomorphism $\Psi_{\tilde{N}(t)}$ acting from
$(tR,\epsilon)\times\Sigma$ into $M$.
If $\tau$ is sufficiently small, $\varphi$
vanishes on $\tilde{N}_u(t)$ and there exists $C(\varphi)>0$, independent of $t$, such that:
\begin{align*}
\big|\,\chi^*\left(\varphi\right)|_{tK_A}\big|&\leq C(\varphi)t;&&\\
\big|\,\chi^*(\varphi)\big(t\Phi_A(r,\sigma)\big)\big|&\leq C(\varphi)tr&&\hspace{-23.5pt}\text{for $r\in\big(R,\textstyle\frac{1}{2}t^{\nu-1}\big)$; and}\\
\left|\,\Psi_{\tilde{N}(t)}^*\left(\varphi|_{\tilde{N}_m(t)}\right)\right|&\leq
C(\varphi)\left(
t^{(1-\nu)(1-\lambda)}+t^{\nu(\mu-1)}\right).&&\\
\intertext{Moreover,}
\big\|\,\varphi|_{\tilde{N}(t)}\big\|_{L^{4/3}}&\leq 
C(\varphi)t^{3\nu}\left(
t^{(1-\nu)(1-\lambda)}+t^{\nu(\mu-1)}\right) && \\
\intertext{and}
\big\|\,\varphi|_{\tilde{N}(t)}\big\|_{L^8_{1,\,0,\,t}}&\leq C(\varphi)\left(t^{(1-\nu)(1-\lambda)}
+t^{\nu(\mu-1)}\right). &&
\end{align*}
\end{prop}

\begin{proof}
The first thing to note is that $\varphi|_{\tilde{N}_u(t)}\equiv0$ as
$\tilde{N}_u(t)$ is a subset of $\hat{N}$ which is coassociative in $(M,\varphi,g)$.
Secondly, $\chi^*(\varphi)=\varphi_0+O(r)$ as $r\rightarrow 0$,
where $\varphi_0$ is given in \eq{s2eq1} and $r$ is the radius function on $B(0;\epsilon_M)$.  Thus, since
$\varphi_0$ vanishes on $A$ as $A$ is coassociative in $\R^7$, we have that
$$\big|\,\chi^*(\varphi)\big(t\Phi_A(t^{-1}r,\sigma)\big)\big|=O(r)$$
for  $r\in (tR,\epsilon)$, so it is certainly true for $r\in
\big(tR,\frac{1}{2}t^{\nu}\big)$. Moreover,
$$\big|\,\chi^*(\varphi)|_{tK_A}\big|=O(t)$$ as $t\rightarrow 0$.

Now recall we have a metric $\tilde{g}(t)=g_{\tilde{N}(t)}$ on $\tilde{N}(t)$, as in the notation of Proposition \ref{s2prop4}, 
and a metric $g_{\text{cone}}$ as given in $\S$\ref{s4subs1}(a). Clearly $\Psi_{\tilde{N}(t)}^*\big(\tilde{g}(t)\big)$ and
$g_{\text{cone}}$ are equivalent on
$(\frac{1}{2}t^{\nu},t^{\nu})\times\Sigma$, under variations in $t$, so any moduli or Levi--Civita connections on this
portion of the cone can be calculated with respect to either metric.  Moreover, the same estimates will hold, perhaps under a
($t$-invariant) rescaling of constants.  We therefore denote the Levi--Civita connection of $g_{\text{cone}}$ by $\nabla$.

We wish to consider $|\Psi_{\tilde{N}(t)}^*(\varphi)|$.  
If we let
$\iota_{\tilde{N}(t)}=\iota:\tilde{C}(t)\rightarrow B(0;\epsilon_M)$ be the inclusion map,
\begin{align*}
|\Psi_{\tilde{N}(t)}^*(\varphi)|&=\big|\,\Psi_{\tilde{N}(t)}^*(\varphi)-\iota_{\tilde{N}(t)}^*\big(\chi^*(\varphi)\big)\big|
+\big|\,\iota_{\tilde{N}(t)}^*\big(\chi^*(\varphi)\big)\big|\\&=
\big|\big(\Phi_{\tilde{N}(t)}-\iota_{\tilde{N}(t)}\big)^*\big(\chi^*(\varphi)\big)\big|+O(r).
\end{align*}
  If
$t$ is sufficiently small so that $t^{\nu}$, the upper bound on
$\big|\chi^*\big(\varphi|_{\tilde{N}_m(t)}\big)\big|$, is small, the maximum of $\big|\Psi_{\tilde{N}(t)}^*\big(\varphi|_{\tilde{N}_m(t)}\big)\big|$
 is determined by
$|\Phi_{\tilde{N}(t)}-\iota_{\tilde{N}(t)}|$ and
$|\nabla(\Phi_{\tilde{N}(t)}-\iota_{\tilde{N}(t)})|$,
 since it acts on the tangent spaces.
Moreover, if the moduli of both maps are small, the dominant terms are \emph{linear} in
the latter.

Recall that
\begin{align*}
\Phi_{\tilde{N}(t)}(r,\sigma)-\iota_{\tilde{N}(t)}(r,\sigma)
&=t\left(1-f_{\text{inc}}(2t^{-\nu}r-1)
\right)\big(\Phi_A(t^{-1}r,\sigma)-\iota(t^{-1}r,\sigma)\big)\\&\quad
+f_{\text{inc}}(2t^{-\nu}r-1)
\big(\Phi_N(r,\sigma)-\iota(r,\sigma)\big).
\end{align*}
Now, by \eq{dt1eq3},
\begin{align*}
\left|t\big(\Phi_A(t^{-1}r,\sigma)-\iota(t^{-1}r,\sigma)\big)\right|&=O\big(t(t^{-1}r)^\lambda\big)=O\big(t^{1-\lambda}
r^\lambda\big)\quad\text{and}\\
\left|t\nabla\big(\Phi_A(t^{-1}r,\sigma)-\iota(t^{-1}r,\sigma)\big)\right|&=
O\big(t\,.\,t^{-1}
(t^{-1}r)^{\lambda-1}\big)=O\big(t^{1-\lambda}r^{\lambda-1}\big)\end{align*}
for $r\in (\frac{1}{2}t^{\nu},t^{\nu})$. Note that there is an extra
$t^{-1}$ factor introduced when taking the derivative. 
For the map $\Phi_N$ we have, by \eq{ch8s1eq1},
$$|\Phi_N(r,\sigma)-\iota(r,\sigma)|=O(r^{\mu})\quad\text{and}\quad \big|\nabla\big(\Phi_N(r,\sigma)-\iota(r,\sigma)\big)\big|=O\big(r^{\mu-1}\big)$$
for $r\in (\frac{1}{2}t^{\nu},t^{\nu})$.

Consider the case where we differentiate the $f_{\text{inc}}$ parts of the expression for $\Phi_{\tilde{N}(t)}-\iota_{\tilde{N}(t)}$.
  These will also produce terms of order $O(t^{1-\lambda}r^{\lambda-1})$ and $O(r^{\mu-1})$ for $r\in (\frac{1}{2}t^{\nu},t^{\nu})$.

We can thus see that, if $\tau$ is sufficiently small,
$$\big|\,\Psi_{\tilde{N}}^*\big(\varphi|_{\tilde{N}_m(t)}\big)\big|=
O\left(\big(t^{-1}r\big)^{\lambda-1}\right)+O\big(r^{\mu-1}\big)$$
for $r\in(\frac{1}{2}t^{\nu},t^{\nu})$.  
  Since $\mu>1$ and $\lambda<1$,
the first set of results follows.

The penultimate result is easily deduced from the fact that the major
contribution to the norm of $\varphi$ is given by its behaviour on
$\tilde{N}_m(t)$ and the fact that $\text{vol}(\tilde{N}_l(t))=O(t^4)$ and $\text{vol}(\tilde{N}_m(t))=O(t^{4\nu})$, in the notation
of Definition \ref{s4dfn2}. 

For the final result we need to repeat all of our estimates for $\big|\nabla_t\varphi|_{\tilde{N}(t)}\big|$, which are very similar
calculations to those we have just detailed.  A key point is that, on $\tilde{N}_m(t)$, its maximum is determined by linear terms
in $|\nabla^2(\Phi_{\tilde{N}(t)}-\iota_{\tilde{N}(t)})|$.  The overall result is that 
$\big|\rho_t\nabla_t\varphi|_{\tilde{N}(t)}\big|$ has the same dependence on $t$ as $\big|\varphi|_{\tilde{N}(t)}\big|$.  
Notice that the dominant terms in our last estimate in the proposition come from the integral over $\tilde{N}_m(t)$, and
$$\int_{\tilde{N}_m(t)}\rho_t^{-4}dV_{\tilde{g}(t)}=O\left(\int_{\frac{1}{2}t^{\nu}}^{t^{\nu}}r^{-4}r^3dr\right)=O(1).$$
Putting together these observations and the earlier calculations completes the proof.
\end{proof}

We can now deduce estimates for $\psi(t)$ from the bounds above.

\begin{prop}\label{s7prop2}  Recall the notation of Definitions \ref{s4dfn2}-\ref{s4dfn4}, \ref{s6dfn5} and \ref{seidfn3} and 
Propositions \ref{seilem2} and \ref{s7prop1}.
The form $\psi(t)$ vanishes on $\tilde{N}_u(t)$, 
\begin{align*}
\|\psi(t)\|_{L^{4/3}_{1,\,-2,\,t}}&\leq C(\varphi)C(\tilde{N},\rho)t^{3\nu}|\log t|^{\frac{1}{2}}
\left(t^{(1-\nu)(1-\lambda)}+t^{\nu(\mu-1)}\right)\\
\intertext{and}
\|\psi(t)\|_{L^8_{2,\,1,\,t}}&\leq C(\varphi)C(\tilde{N},\rho)\big(1+t^{3(\nu-1)}\big)
\left(t^{(1-\nu)(1-\lambda)}+t^{\nu(\mu-1)}\right)
\end{align*}
\end{prop}

\begin{proof}
Since $\psi(t)\in V^8_{2}(t)$ by Corollary \ref{seicor1}, using the notation there, 
it satisfies the conditions of Proposition \ref{seilem2}.  Recall that $d\big(\psi(t)\big)=\varphi|_{\tilde{N}(t)}$.
Therefore, 
\begin{align*}
\|\psi(t)\|_{L^{4/3}_{1,\,-2,\,t}}
&\leq C(\tilde{N},\rho)|\log t|^{\frac{1}{2}}\big\|\,\varphi|_{\tilde{N}(t)}\big\|_{L^{4/3}}\\
\intertext{and}
\|\psi(t)\|_{L^{8}_{2,\,1,\,t}}
&\leq C(\tilde{N},\rho)\left(\big\|\,\varphi|_{\tilde{N}(t)}\big\|_{L^8_{1,\,0,\,t}}+
t^{-3}\big\|\,\varphi|_{\tilde{N}(t)}\big\|_{L^{4/3}}\right).
\end{align*}
Using Proposition \ref{s7prop1} gives the result.
\end{proof}

\noindent Following on from this proposition, we prove a lemma which will be useful later.

\begin{lemma}\label{s7lem1} Use the notation of Definitions \ref{s4dfn2} and \ref{s5dfn3} 
and Propositions \ref{seilem2} and \ref{s7prop1}.
If $\tau$ is sufficiently small then, for all $t\in(0,\tau)$,
\begin{align*}
C(\varphi)C(\tilde{N},\rho)t^{3\nu}|\log t|^{\frac{1}{2}}
\left(
t^{(1-\nu)(1-\lambda)}+t^{\nu(\mu-1)}\right)&\leq \frac{\tilde{\epsilon}t^3}{2}\\
\intertext{and}
C(\varphi)C(\tilde{N},\rho)\big(1+t^{3(\nu-1)}\big)
\left(t^{(1-\nu)(1-\lambda)}+t^{\nu(\mu-1)}\right)&\leq \frac{\tilde{\epsilon}}{2}\,.
\end{align*}
\end{lemma}

\begin{proof}  
Since $\nu\in(0,1)$, $\lambda<1$ and $\mu>1$, it is clear that $(1-\nu)(1-\lambda)>0$ and $\nu(\mu-1)>0$. Further, recall
$t^\kappa\log t\rightarrow 0$ as $t\rightarrow 0$ for any $\kappa>0$.  It is therefore enough to
show that $3(\nu-1)+(1-\nu)(1-\lambda)>0$ and $3(\nu-1)+\nu(\mu-1)>0$ to prove the result for sufficiently small $\tau$.
First,
$$3(\nu-1)+(1-\nu)(1-\lambda)=(\nu-1)(\lambda+2)>0$$
since $\lambda<-2$ and $\nu<1$.  Second,
$$3(\nu-1)+\nu(\mu-1)=\nu(\mu+2)-3>0$$
by the choice of $\nu>\frac{3}{\mu+2}$.
\end{proof}

\subsection{Estimates involving {\boldmath $P_{H_t}$}}

We begin with a result which is elementary given our earlier work.  Recall that if $\|\alpha\|_{L^8_{2,\,1,\,t}}\leq\tilde{\epsilon}$,
then we can ensure that $\|\alpha\|_{C^1_{1,\,t}}$ is small by making $\tilde{\epsilon}$ sufficiently small.

\begin{prop}\label{s7prop6}  Use the notation of Propositions \ref{s6prop1}, \ref{s6prop3}, \ref{s6prop2} and \ref{seilem2}, Corollary \ref{seicor1}
 and Definitions \ref{s5dfn3} and \ref{seidfn3}.
 If $\tau$ is sufficiently small and $\alpha,\beta\in V^8_{2}(t)$ with $\|\alpha\|_{L^8_{2,\,1,\,t}},
\|\beta\|_{L^8_{2,\,1,\,t}}\leq \tilde{\epsilon}$,
\begin{align*}
\|P_{H_t}(\alpha)&\|_{L^{4/3}_{1,\,-2,\,t}}\leq C(\tilde{N},\rho)p_0|\log t|^{\frac{1}{2}}
\|\alpha\|_{L^{4/3}_{1,\,-2,\,t}}^2,\\
\|P_{H_t}(\alpha)&\|_{L^{8}_{2,\,1,\,t}}\leq C(\tilde{N},\rho)\left(p_1\|\alpha\|_{L^8_{2,\,1,\,t}}^2+
p_0t^{-3}\|\alpha\|_{L^{4/3}_{1,\,-2,\,t}}^2\right),\\
\|P_{H_t}(\alpha)&-P_{H_t}(\beta)\|_{L^{4/3}_{1,\,-2,\,t}}\\
&\leq C(\tilde{N},\rho)C(P_{F_t})|\log t|^{\frac{1}{2}}\|\alpha-\beta\|_{L^{4/3}_{1,\,-2,\,t}}
\big(\|\alpha\|_{L^{4/3}_{1,\,-2,\,t}}+\|\beta\|_{L^{4/3}_{1,\,-2,\,t}}\big)\;\,\text{and}\\
\|P_{H_t}(\alpha)&-P_{H_t}(\beta)\|_{L^8_{2,\,1,\,t}}\\
&\leq C(\tilde{N},\rho)C(P_{F_t})\bigg(\|\alpha-\beta\|_{L^8_{2,\,1,\,t}}
\big(\|\alpha\|_{L^8_{2,\,1,\,t}}+\|\beta\|_{L^8_{2,\,1,\,t}}\big)\\&\qquad\qquad\qquad\qquad+
t^{-3}\|\alpha-\beta\|_{L^{4/3}_{1,\,-2,\,t}}
\big(\|\alpha\|_{L^{4/3}_{1,\,-2,\,t}}+\|\beta\|_{L^{4/3}_{1,\,-2,\,t}}\big)\bigg).
\end{align*}
\end{prop}

\begin{proof}
The first two inequalities are proved by applying Propositions \ref{seilem2} and \ref{s6prop1}.  The last two inequalities follow 
again from Proposition \ref{seilem2},
but then one applies Proposition \ref{s6prop3}.
\end{proof}

\noindent We follow our estimate with a useful, though easy, lemma.

\begin{lemma}\label{s7lem2}  Use the notation of Definition \ref{s5dfn3} and Propositions \ref{s6prop1}, \ref{s6prop3} and \ref{seilem2}.
If $\tau$ and $\tilde{\epsilon}$ are sufficiently small then, for all $t\in(0,\tau)$,
\begin{gather*}
C(\tilde{N},\rho)p_0|\log t|^{\frac{1}{2}}\tilde{\epsilon}^2t^6\leq\frac{\tilde{\epsilon}t^3}{2},
\\
C(\tilde{N},\rho)\tilde{\epsilon}^2\left(p_1+
p_0t^{3}\right)\leq\frac{\tilde{\epsilon}}{2}\;\,\text{and}\\
0<2C(\tilde{N},\rho)C(P_{F_t})\tilde{\epsilon}\left(1+|\log t|^{\frac{1}{2}}t^3\right)< 1
\end{gather*}
\end{lemma}

\begin{proof}
This result follows from the observation that $t^3|\log t|^{\frac{1}{2}}\rightarrow 0$ as $t\rightarrow 0$.
\end{proof}

\subsection{The contraction map}

The idea is that we want a form $\alpha\in
L^8_{2}\big(\tilde{U}(t)\big)$ which satisfies $F_t(\alpha)=0$.  We
know, however, that it is enough to find $\alpha$ satisfying
$H_t(\alpha)=0$ so that its graph defines a coassociative
(nonsingular) deformation of $\tilde{N}(t)$.  We start with a definition.

\begin{dfn}\label{s7dfn2} 
In the notation of Definition \ref{s5dfn3} and Corollary \ref{seicor1}, define
$$X(t)=\{\alpha\in V^8_{2}(t)\,:\,\|\alpha\|_{L^{4/3}_{1,\,-2,\,t}}\leq \tilde{\epsilon}t^3\;\,\text{and}
\;\,\|\alpha\|_{L^8_{2,\,1,\,t}}\leq \tilde{\epsilon}\}$$
and endow it with the norm
$$\|\alpha\|_{X(t)}=\|\alpha\|_{L^{4/3}_{1,\,-2,\,t}}+\|\alpha\|_{L^8_{2,\,1,\,t}}.$$
Notice that it is a \emph{closed} subset of a Banach space and
therefore \emph{complete}, and that its elements have their graphs
contained in $\tilde{U}(t)$ by Definition \ref{s5dfn3}.  
\end{dfn}

With this definition, if $\alpha\in X(t)$, we can write the map $H_t$ in the
form
$$H_t(\alpha)=\psi(t)+\alpha+P_{H_t}(\alpha),$$
noting that $\alpha_{H_t}=\alpha$, in the notation of Definition \ref{seidfn3}, since $\alpha$ lies in $V^8_{2}(t)$.
Hence, if we define
\begin{equation}\label{ct}
\mathcal{C}_t(\alpha)=\psi(t)-P_{H_t}(\alpha),
\end{equation}
$\alpha$ must satisfy $\mathcal{C}_t(\alpha)=\alpha$.  Clearly $\mathcal{C}_t$ maps $X(t)$ into $V^8_2(t)$.  
We prove that $\mathcal{C}_t$ is in fact a \emph{contraction}.  Applying the
Contraction Mapping Theorem will give us a form in $X(t)$ satisfying
$\mathcal{C}_t(\alpha)=\alpha$ and thus $H_t(\alpha)=0$.


We first need to show that if $\alpha\in X(t)$
then 
$\mathcal{C}_t(\alpha)\in X(t)$.  
Choose $\tau>0$ and $\tilde{\epsilon}>0$ sufficiently small. 
One may quickly see that, by Propositions \ref{s7prop2} and \ref{s7prop6} and Lemmas \ref{s7lem1} and \ref{s7lem2}, if 
$\alpha\in X(t)$,
\begin{align*}
&\|\mathcal{C}_t(\alpha)\|_{L^{4/3}_{1,\,-2,\,t}}\leq 
\|\psi(t)\|_{L^{4/3}_{1,\,-2,\,t}}+\|P_{H_t}(\alpha)\|_{L^{4/3}_{1,\,-2,\,t}}\\
&\leq 
C(\varphi)C(\tilde{N},\rho)|\log t|^{\frac{1}{2}}t^{3\nu}
\left(
t^{(1-\nu)(1-\lambda)}+t^{\nu(\mu-1)}\right)+C(\tilde{N},\rho)p_0|\log t|^{\frac{1}{2}}\tilde{\epsilon}^2t^6\\
&\leq \tilde{\epsilon}t^3.
\end{align*}
Using the same results we can show that $\|\mathcal{C}_t(\alpha)\|_{L^8_{2,\,1,\,t}}\leq\tilde{\epsilon}$.

With this observation, we have our first necessary result. 

\begin{prop}\label{s7prop3} In the notation of Definitions \ref{s5dfn3} and \ref{s7dfn2}, the map 
$\mathcal{C}_t$ given in \eq{ct} takes $X(t)$ to
itself for all $t\in(0,\tau)$, if $\tau$ and $\tilde{\epsilon}$ are sufficiently small.
\end{prop}

\noindent We now turn our attention to a more important proposition.

\begin{prop}\label{s7prop4}
In the notation of Definitions \ref{s5dfn3} and \ref{s7dfn2}, if $\tau$ and $\tilde{\epsilon}>0$ are sufficiently small,
$\mathcal{C}_t:X(t)\rightarrow X(t)$, given in \eq{ct}, is a contraction for all $t\in(0,\tau)$.
\end{prop}

\begin{proof}
Let $\alpha,\beta\in X(t)$.  Then
$$\mathcal{C}_t(\alpha)-\mathcal{C}_t(\beta)=P_{H_t}(\beta)-P_{H_t}(\alpha),$$
where $P_{H_t}$ is given by Definition \ref{seidfn3}.
We also know that, by Proposition \ref{s7prop6},
\begin{align*}
\|P_{H_t}(\alpha)&-P_{H_t}(\beta)\|_{L^{4/3}_{1,\,-2,\,t}} \\
&\leq C(\tilde{N},\rho)C(P_{F_t})|\log t|^{\frac{1}{2}}\|\alpha-\beta\|_{L^{4/3}_{1,\,-2,\,t}}
\big(\|\alpha\|_{L^{4/3}_{1,\,-2,\,t}}+\|\beta\|_{L^{4/3}_{1,\,-2,\,t}}\big)\\
&\leq 2C(\tilde{N},\rho)C(P_{F_t})\tilde{\epsilon}|\log t|^{\frac{1}{2}}t^3\|\alpha-\beta\|_{L^{4/3}_{1,\,-2,\,t}},
\end{align*}
where $C(P_{F_t})$ is given by Proposition \ref{s6prop3} and $C(\tilde{N},\rho)$ is given by Proposition \ref{seilem2}.
Using Proposition \ref{s7prop6} again,
\begin{align*}
\|P_{H_t}(\alpha)&-P_{H_t}(\beta)\|_{L^{8}_{2,\,1,\,t}} \\
&\leq C(\tilde{N},\rho)C(P_{F_t})\bigg(\|\alpha-\beta\|_{L^8_{2,\,1,\,t}}
\big(\|\alpha\|_{L^8_{2,\,1,\,t}}+\|\beta\|_{L^8_{2,\,1,\,t}}\big)\\&\qquad\qquad\qquad\qquad+
t^{-3}\|\alpha-\beta\|_{L^{4/3}_{1,\,-2,\,t}}
\big(\|\alpha\|_{L^{4/3}_{1,\,-2,\,t}}+\|\beta\|_{L^{4/3}_{1,\,-2,\,t}}\big)\bigg)\\
&\leq 2C(\tilde{N},\rho)C(P_{F_t})\tilde{\epsilon}\|\alpha-\beta\|_{X(t)}.
\end{align*}
Thus,
$$\|P_{H_t}(\alpha)-P_{H_t}(\beta)\|_{X(t)}\leq 2C(\tilde{N},\rho)C(P_{F_t})\tilde{\epsilon}
\big(1+|\log t|^{\frac{1}{2}}t^3\big)\|\alpha-\beta\|_{X(t)}.$$
Applying Lemma \ref{s7lem2}, there exists $\kappa(\mathcal{C}_t)\in(0,1)$, independent of $\alpha$ and $\beta$, such that
$$\|\mathcal{C}_t(\alpha)-\mathcal{C}_t(\beta)\|_{X(t)}\leq\kappa(\mathcal{C}_t)\|\alpha-\beta\|_{X(t)}.$$
\end{proof}

The work of this section culminates in our final proposition.

\begin{prop}\label{s7prop5}  Use the notation of Definitions \ref{s5dfn3}, \ref{s6dfn1}, \ref{seidfn3} and \ref{s7dfn2}, 
and suppose $\tau$ and $\tilde{\epsilon}$ are sufficiently small.
For all $t\in(0,\tau)$, there exists a
(unique) solution $\alpha(t)$ in 
$X(t)$ to $H_t\big(\alpha(t)\big)=0$, which thus solves
$F_t\big(\alpha(t)\big)=0$.
\end{prop}

\begin{proof}
By Proposition
\ref{s7prop4}, the map $\mathcal{C}_t$ from
$X(t)$ to itself
is a contraction.  Applying the Contraction Mapping Theorem,
recalling that $\mathcal{C}_t$ acts on a complete space, there exists
a (unique) $\alpha(t)\in X(t)$ satisfying
$$\alpha(t)=\mathcal{C}_t\big(\alpha(t)\big)=-\psi(t)-P_{H_t}\big(\alpha(t)\big)=-H_t\big(\alpha(t)\big)+\alpha(t)$$
by Proposition \ref{s6prop2}.
Therefore, $H_t\big(\alpha(t)\big)=0$, which implies $F_t\big(\alpha(t)\big)=0$.
\end{proof}

\section{Desingularization: final stage}\label{s9}

Our final section deals quickly
with the \emph{regularity} of $\alpha(t)$ and gives the principal result of the paper.

\begin{prop}\label{s9prop1}  In the notation of Propositions \ref{s5prop3} and \ref{s7prop5}, $\alpha(t)$ lies in $C^{\infty}\big(\tilde{U}(t)\big)$.
\end{prop}

\begin{proof} Recall the definition of $\tilde{\epsilon}>0$ in Definition \ref{s5dfn3} and that the form $\alpha(t)\in 
L^8_2\big((\Lambda^2_+)_{\tilde{g}(t)}T^*\tilde{N}(t)\big)$ with $\|\alpha(t)\|_{L^8_{2,\,1,\,t}}\leq\tilde{\epsilon}$.
Therefore, by Definition \ref{s5dfn3}, $\alpha(t)\in L^8_2\big(\tilde{U}(t)\big)$ and we can ensure $\|\alpha\|_{C^1_{1,\,t}}$ is small by making $\tilde{\epsilon}$ smaller. 
Since $F_t\big(\alpha(t)\big)=0$, we can apply Proposition \ref{s8prop2} 
to give the result.  
\end{proof}

The work in this article allows us to conclude with our key theorem.

\begin{thm}\label{s9thm1}  Use the notation of $\S$\ref{s4subs1}(a)-(c) and Definitions \ref{s4dfn2} and \ref{TN}.  
In particular, recall we have a coassociative 4-fold $N$ in a $\varphi$-closed 7-manifold $(M,\varphi,g)$ with a singularity modelled 
on a cone $C$, with link $\Sigma$ such that $b^1(\Sigma)=0$, and a coassociative 4-fold $A$ in $\R^7$ asymptotically conical to $C$
with rate $\lambda<-2$.  Further, we have, for each $t\in(0,\tau)$ where $\tau>0$ is sufficiently small, 
a nonsingular compact 4-fold $\tilde{N}(t)$ in $M$ which is formed by 
gluing $tA$ into $N$ at the singularity.

For all $t\in(0,\tau)$, there exists a nonsingular, compact, coassociative
deformation $N^{\prime}(t)$ of $\tilde{N}(t)$, which lies in $T_N$, an open neighbourhood of $N$ in $M$.  
Moreover, the family $\{N^{\prime}(t)\,:\,t\in(0,\tau)\}$
is smooth in $t$ and $N^{\prime}(t)$ converges in the sense of currents to $N$ as $t\rightarrow 0$, as described in Definition \ref{s1dfn1}.
\end{thm}

\begin{proof} Use the notation of Proposition \ref{s5prop3} and Definition \ref{s6dfn1}. 
By Propositions \ref{s7prop5} and \ref{s9prop1}, for all $t\in(0,\tau)$, we have $\alpha(t)\in C^{\infty}\big(\tilde{U}(t)\big)$ satisfying $F_t\big(\alpha(t)\big)=0$.  
Define 
$$N^{\prime}(t)=\tilde{N}_{\alpha(t)}(t)=\big(\tilde{\delta}(t)\circ\pi_{\alpha(t)}\big)\big(G_{\alpha(t)}\big),$$
where $G_{\alpha(t)}$ is the graph of $\alpha(t)$.  We see that $N^{\prime}(t)$ 
is a compact coassociative deformation of $\tilde{N}(t)$ which 
lies in $T_N$ by construction.   

By the note made after Definition \ref{s4dfn2}, the family of $\tilde{N}(t)$ is smooth in $t$ and converges to
$N$ in the sense of currents as $t\rightarrow 0$.  
Clearly the family of $N^{\prime}(t)$ is thus also smooth in $t$.  Recall that $\alpha(t)\in X(t)$, given in Definition 
\ref{s7dfn2}, and hence has bounds imposed on its norm in certain weighted Sobolev spaces.  These bounds ensure that, 
for any compactly supported 4-form $\xi$ on $M$, $$\int_{\tilde{N}(t)}\xi-\int_{N^{\prime}(t)}\xi\rightarrow 0\quad\text{as 
$t\rightarrow 0$.}$$
Therefore $N^{\prime}(t)\rightarrow N$ as $t\rightarrow 0$ in the sense of currents, described in Definition
\ref{s1dfn1}, as claimed.
\end{proof}

\noindent This theorem gives us our smooth family of coassociative desingularizations of
$N$, the topology and geometry of which can be explicitly described by considering the construction of $\tilde{N}(t)$ 
described in Definition \ref{s4dfn2}.

\begin{remark}
All of the work in this article can easily be extended to CS coassociative 4-folds $N$ with $s$ singularities, 
by gluing in $s$ AC coassociative 4-folds $A_1,\ldots, A_s$ in $\R^7$.  
It is clear that the analysis will be almost the same as for
the single singularity scenario we have studied here.  
\end{remark}

\noindent By this remark, we realise that Definition \ref{s4dfn2} and Theorem \ref{s9thm1} imply our main result Theorem \ref{s1thm1}
as intended. 

\bigskip

\noindent{\bf Acknowledgements}\quad The author would like to thank Dominic Joyce and Alexei Kovalev for their helpful comments and
advice.  In particular, to Alexei Kovalev, for bringing the references \cite{Kovalev2} and \cite{Mazya} to the author's attention.  
The author is also extremely grateful to the referee for providing a detailed and useful report.

\end{document}